\documentclass[12pt]{article}

\usepackage{amsmath,amsfonts,amssymb,graphicx,a4,bm,subeqnarray}
\usepackage{indent}
\usepackage{color}
\usepackage{url}       

\def\theequation{\thesection.\arabic{equation}}
\newcommand{\zeq}{\setcounter{equation}{0}}

\definecolor{GreenYellow}{cmyk}{0.15,0,0.69,0}
\definecolor{Yellow}{cmyk}{0,0,1,0}
\definecolor{Goldenrod}{cmyk}{0,0.10,0.84,0}
\definecolor{Dandelion}{cmyk}{0,0.29,0.84,0}
\definecolor{Apricot}{cmyk}{0,0.32,0.52,0}
\definecolor{Peach}{cmyk}{0,0.50,0.70,0}
\definecolor{Melon}{cmyk}{0,0.46,0.50,0}
\definecolor{YellowOrange}{cmyk}{0,0.42,1,0}
\definecolor{Orange}{cmyk}{0,0.61,0.87,0}
\definecolor{BurntOrange}{cmyk}{0,0.51,1,0}
\definecolor{Bittersweet}{cmyk}{0,0.75,1,0.24}
\definecolor{RedOrange}{cmyk}{0,0.77,0.87,0}
\definecolor{Mahogany}{cmyk}{0,0.85,0.87,0.35}
\definecolor{Maroon}{cmyk}{0,0.87,0.68,0.32}
\definecolor{BrickRed}{cmyk}{0,0.89,0.94,0.28}
\definecolor{Red}{cmyk}{0,1,1,0}
\definecolor{OrangeRed}{cmyk}{0,1,0.50,0}
\definecolor{RubineRed}{cmyk}{0,1,0.13,0}
\definecolor{WildStrawberry}{cmyk}{0,0.96,0.39,0}
\definecolor{Salmon}{cmyk}{0,0.53,0.38,0}
\definecolor{CarnationPink}{cmyk}{0,0.63,0,0}
\definecolor{Magenta}{cmyk}{0,1,0,0}
\definecolor{VioletRed}{cmyk}{0,0.81,0,0}
\definecolor{Rhodamine}{cmyk}{0,0.82,0,0}
\definecolor{Mulberry}{cmyk}{0.34,0.90,0,0.02}
\definecolor{RedViolet}{cmyk}{0.07,0.90,0,0.34}
\definecolor{Fuchsia}{cmyk}{0.47,0.91,0,0.08}
\definecolor{Lavender}{cmyk}{0,0.48,0,0}
\definecolor{Thistle}{cmyk}{0.12,0.59,0,0}
\definecolor{Orchid}{cmyk}{0.32,0.64,0,0}
\definecolor{DarkOrchid}{cmyk}{0.40,0.80,0.20,0}
\definecolor{Purple}{cmyk}{0.45,0.86,0,0}
\definecolor{Plum}{cmyk}{0.50,1,0,0}
\definecolor{Violet}{cmyk}{0.79,0.88,0,0}
\definecolor{RoyalPurple}{cmyk}{0.75,0.90,0,0}
\definecolor{BlueViolet}{cmyk}{0.86,0.91,0,0.04}
\definecolor{Periwinkle}{cmyk}{0.57,0.55,0,0}
\definecolor{CadetBlue}{cmyk}{0.62,0.57,0.23,0}
\definecolor{CornflowerBlue}{cmyk}{0.65,0.13,0,0}
\definecolor{MidnightBlue}{cmyk}{0.98,0.13,0,0.43}
\definecolor{NavyBlue}{cmyk}{0.94,0.54,0,0}
\definecolor{RoyalBlue}{cmyk}{1,0.50,0,0}
\definecolor{Blue}{cmyk}{1,1,0,0}
\definecolor{Cerulean}{cmyk}{0.94,0.11,0,0}
\definecolor{Cyan}{cmyk}{1,0,0,0}
\definecolor{ProcessBlue}{cmyk}{0.96,0,0,0}
\definecolor{SkyBlue}{cmyk}{0.62,0,0.12,0}
\definecolor{Turquoise}{cmyk}{0.85,0,0.20,0}
\definecolor{TealBlue}{cmyk}{0.86,0,0.34,0.02}
\definecolor{Aquamarine}{cmyk}{0.82,0,0.30,0}
\definecolor{BlueGreen}{cmyk}{0.85,0,0.33,0}
\definecolor{Emerald}{cmyk}{1,0,0.50,0}
\definecolor{JungleGreen}{cmyk}{0.99,0,0.52,0}
\definecolor{SeaGreen}{cmyk}{0.69,0,0.50,0}
\definecolor{Green}{cmyk}{1,0,1,0}
\definecolor{ForestGreen}{cmyk}{0.91,0,0.88,0.12}
\definecolor{PineGreen}{cmyk}{0.92,0,0.59,0.25}
\definecolor{LimeGreen}{cmyk}{0.50,0,1,0}
\definecolor{YellowGreen}{cmyk}{0.44,0,0.74,0}
\definecolor{SpringGreen}{cmyk}{0.26,0,0.76,0}
\definecolor{OliveGreen}{cmyk}{0.64,0,0.95,0.40}
\definecolor{RawSienna}{cmyk}{0,0.72,1,0.45}
\definecolor{Sepia}{cmyk}{0,0.83,1,0.70}
\definecolor{Brown}{cmyk}{0,0.81,1,0.60}
\definecolor{Tan}{cmyk}{0.14,0.42,0.56,0}
\definecolor{Gray}{cmyk}{0,0,0,0.50}
\definecolor{Black}{cmyk}{0,0,0,1}
\definecolor{White}{cmyk}{0,0,0,0}

\begin{document}

\def\blu{\color{Blue}}
\def\mag{\color{Maroon}}
\def\red{\color{Red}}
\def\green{\color{ForestGreen}}
\def\orange{\color{RedOrange}}



\def\\{\noindent}
\def\PP{{\mathcal P}}
\def\I{{\rm I}}
\def\0{\emptyset}
\def\vv{\vskip.15cm}

\def\La{\Lambda}

%
%
\newcommand{\stirlingsubset}[2]{\genfrac{\{}{\}}{0pt}{}{#1}{#2}}
\newcommand{\stirlingcycle}[2]{\genfrac{[}{]}{0pt}{}{#1}{#2}}
\newcommand{\assocstirlingsubset}[3]{%
{\genfrac{\{}{\}}{0pt}{}{#1}{#2}}_{\! \ge #3}}
\newcommand{\assocstirlingcycle}[3]{{\genfrac{[}{]}{0pt}{}{#1}{#2}}_{\ge #3}}
\newcommand{\euler}[2]{\genfrac{\langle}{\rangle}{0pt}{}{#1}{#2}}
\newcommand{\eulergen}[3]{{\genfrac{\langle}{\rangle}{0pt}{}{#1}{#2}}_{\! #3}}
\newcommand{\eulersecond}[2]{\left\langle\!\! \euler{#1}{#2} \!\!\right\rangle}
\newcommand{\eulersecondgen}[3]{%
{\left\langle\!\! \euler{#1}{#2} \!\!\right\rangle}_{\! #3}}
\newcommand{\binomvert}[2]{\genfrac{\vert}{\vert}{0pt}{}{#1}{#2}}

\title{\vspace*{-2.5cm} Complex zero-free regions at large $|q|$ \\
        for multivariate Tutte polynomials \\
        (alias Potts-model partition functions) \\
        with general complex edge weights}

\author{
  \\[-8mm]
  {\small Bill Jackson}                                    \\[-2mm]
  {\small\it School of Mathematical Sciences}  \\[-2mm]
  {\small\it Queen Mary University of London} \\[-2mm]
  {\small\it Mile End Road} \\[-2mm]
  {\small\it London E1 4NS, England}                         \\[-2mm]
  {\small\tt B.JACKSON@QMUL.AC.UK}                        \\[4mm]
{\small Aldo Procacci}                                    \\[-2mm]
  {\small\it Departamento de Matem\'atica}  \\[-2mm]
  {\small\it Universidade Federal de Minas Gerais} \\[-2mm]
  {\small\it Av. Ant\^onio Carlos, 6627 -- Caixa Postal 702 } \\[-2mm]
  {\small\it 30161-970 Belo Horizonte, MG -- BRASIL}           \\[-2mm]
  {\small\tt ALDO@MAT.UFMG.BR}\\[4mm]
{\small Alan D.~Sokal\thanks{Also at Department of Mathematics,
           University College London, London WC1E 6BT, England.}}   \\[-2mm]
  {\small\it Department of Physics}       \\[-2mm]
  {\small\it New York University}         \\[-2mm]
  {\small\it 4 Washington Place}          \\[-2mm]
  {\small\it New York, NY 10003 USA}      \\[-2mm]
  {\small\tt SOKAL@NYU.EDU}               \\[-7mm]
  {\protect\makebox[5in]{\quad}}  
}

\date{Version 1: October 26, 2008 \\
      Version 2: November 20, 2009 \\[1.5mm]
      Version 3: September 23, 2011
     }
\maketitle

\vspace*{-8mm}

\begin{abstract}
We find zero-free regions in the complex plane at large $|q|$
for the multivariate Tutte polynomial
(also known in statistical mechanics as the Potts-model partition function)
$Z_{G}(q,\bm w)$ of a graph $G$
with general complex edge weights $\bm w = \{w_e\}$.
This generalizes a result of Sokal \cite{Sokal_01}
that applies only within the complex antiferromagnetic regime $|1+w_e| \le 1$.
Our proof uses the polymer-gas representation of the multivariate Tutte
polynomial together with the Penrose identity.
\end{abstract}


\bigskip
\noindent
{\bf Mathematics Subject Classification (MSC) codes:}
05C15 (Primary); 05A20, 05B35, 05C99, 05E99, 30C15, 82B20 (Secondary).

\clearpage

\newtheorem{defin}{Definition}[section]
\newtheorem{definition}{Definition}[section]
\newtheorem{prop}[defin]{Proposition}
\newtheorem{proposition}[defin]{Proposition}
\newtheorem{lem}[defin]{Lemma}
\newtheorem{lemma}[defin]{Lemma}
\newtheorem{guess}[defin]{Conjecture}
\newtheorem{ques}[defin]{Question}
\newtheorem{question}[defin]{Question}
\newtheorem{prob}[defin]{Problem}
\newtheorem{problem}[defin]{Problem}
\newtheorem{thm}[defin]{Theorem}
\newtheorem{theorem}[defin]{Theorem}
\newtheorem{cor}[defin]{Corollary}
\newtheorem{corollary}[defin]{Corollary}
\newtheorem{conj}[defin]{Conjecture}
\newtheorem{conjecture}[defin]{Conjecture}
\newtheorem{examp}[defin]{Example}
\newtheorem{example}[defin]{Example}

\newtheorem{pro}{Problem}
\newtheorem{clm}{Claim}
\newtheorem{con}{Conjecture}

\renewcommand{\theenumi}{\alph{enumi}}
\renewcommand{\labelenumi}{(\theenumi)}
\def\prf{\par\noindent{\bf Proof.\enspace}\rm}
\def\rmk{\par\medskip\noindent{\bf Remark.\enspace}\rm}

\newcommand{\be}{\begin{equation}}
\newcommand{\ee}{\end{equation}}
\newcommand{\<}{\langle}
\renewcommand{\>}{\rangle}
\newcommand{\widebar}{\overline}
\def\reff#1{(\protect\ref{#1})}
\def\spose#1{\hbox to 0pt{#1\hss}}
\def\ltapprox{\mathrel{\spose{\lower 3pt\hbox{$\mathchar"218$}}
 \raise 2.0pt\hbox{$\mathchar"13C$}}}
\def\gtapprox{\mathrel{\spose{\lower 3pt\hbox{$\mathchar"218$}}
 \raise 2.0pt\hbox{$\mathchar"13E$}}}
\def\textprime{${}^\prime$}
\def\proof{\par\medskip\noindent{\sc Proof.\ }}
\def\sketchofproof{\par\medskip\noindent{\sc Sketch of Proof.\ }}
\newcommand{\qed}{\quad $\Box$ \medskip \medskip}
\def\proofof#1{\bigskip\noindent{\sc Proof of #1.\ }}
\def\altproofof#1{\bigskip\noindent{\sc Alternate proof of #1.\ }}
\def\half{ {1 \over 2} }
\def\third{ {1 \over 3} }
\def\twothird{ {2 \over 3} }
\def\smfrac#1#2{\textstyle{#1\over #2}}
\def\smhalf{ \smfrac{1}{2} }
\newcommand{\real}{\mathop{\rm Re}\nolimits}
\renewcommand{\Re}{\mathop{\rm Re}\nolimits}
\newcommand{\imag}{\mathop{\rm Im}\nolimits}
\renewcommand{\Im}{\mathop{\rm Im}\nolimits}
\newcommand{\sgn}{\mathop{\rm sgn}\nolimits}
\def\hboxscript#1{ {\hbox{\scriptsize\em #1}} }

\newcommand{\restrict}{\upharpoonright}
\renewcommand{\emptyset}{\varnothing}

\newcommand{\scra}{{\mathcal{A}}}
\newcommand{\scrb}{{\mathcal{B}}}
\newcommand{\scrc}{{\mathcal{C}}}
\newcommand{\scrf}{{\mathcal{F}}}
\newcommand{\scrg}{{\mathcal{G}}}
\newcommand{\scrh}{{\mathcal{H}}}
\newcommand{\scrk}{{\mathcal{K}}}
\newcommand{\scrl}{{\mathcal{L}}}
\newcommand{\scrm}{{\mathcal{M}}}
\newcommand{\scro}{{\mathcal{O}}}
\newcommand{\scrp}{{\mathcal{P}}}
\newcommand{\scrr}{{\mathcal{R}}}
\newcommand{\scrs}{{\mathcal{S}}}
\newcommand{\scrt}{{\mathcal{T}}}
\newcommand{\scrv}{{\mathcal{V}}}
\newcommand{\scrw}{{\mathcal{W}}}
\newcommand{\scrz}{{\mathcal{Z}}}

\newcommand{\sV}{{\mathsf{V}}}
\newcommand{\sE}{{\mathsf{E}}}

\newcommand{\w}{{\bm{w}}}
\newcommand{\wtilde}{{\widetilde{\bm{w}}}}
\newcommand{\what}{{\widehat{\bm{w}}}}
\newcommand{\wdoubleprimex}{w^{[x]}}
\newcommand{\bwdoubleprimex}{{\bm{w}^{[x]}}}
\newcommand{\wtildex}{{\widetilde{w}^{[x]}}}
\newcommand{\bwtildex}{{\widetilde{\bm{w}}^{[x]}}}
\newcommand{\Rtilde}{{\widetilde{\bf R}}}
\newcommand{\Rhat}{{\widehat{\bf R}}}
\newcommand{\smalln}{\mbox{\scriptsize\bf n}}
\newcommand{\blambda}{\mbox{\boldmath $\lambda$}}
\newcommand{\smallblambda}{\mbox{\scriptsize\boldmath $\lambda$}}
\newcommand{\btheta}{\mbox{\boldmath $\theta$}}
\newcommand{\balpha}{\mbox{\boldmath $\alpha$}}
\newcommand{\bdelta}{\mbox{\boldmath $\delta$}}
\newcommand{\smallbalpha}{\mbox{\scriptsize\boldmath $\alpha$}}
\def\twotilde#1{{\widetilde{\widetilde{#1}}}}
\def\twohat#1{{\widehat{\widehat{#1}}}}

\newcommand{\C}{{\mathbb C}}
\newcommand{\Z}{{\mathbb Z}}
\newcommand{\N}{{\mathbb N}}
\newcommand{\R}{{\mathbb R}}
\newcommand{\Cbar}{{\overline{\C}}}

\newcommand{\Zhat}{\widehat{Z}}
\newcommand{\Ztilde}{\widetilde{Z}}
\newcommand{\Ptilde}{\widetilde{P}}
\newcommand{\Ctilde}{\widetilde{C}}
\newcommand{\series}{{\,\bowtie_q\,}}
\newcommand{\seriesq}{{\,\bowtie_q\,}}
\newcommand{\seriesnoq}{{\,\bowtie\,}}
\renewcommand{\parallel}{\Vert}

%
%
\font\fourrm  = cmr5 
\def\sspose#1{\hbox to 0pt{#1\hss}}
\def\gtalmost{\mathrel{\sspose{\lower 0.75pt\hbox{\kern-1.5pt\fourrm (\quad)}}
  \raise 2.0pt\hbox{$\ge$}}}
\def\ltalmost{\mathrel{\sspose{\lower 0.75pt\hbox{\kern-1.5pt\fourrm (\quad)}}
  \raise 2.0pt\hbox{$\le$}}}

{\centerline{\Large}}

\section{Introduction}
\zeq

A decade ago, Sokal \cite{Sokal_01} proved that
if $G=(V,E)$ is a loopless graph\footnote{
   All graphs in this paper are finite and undirected;
   furthermore, they are {\em allowed}\/ to contain loops and multiple edges
   unless we explicitly state otherwise.
}
of maximum degree $\Delta$,
then all the roots (real or complex) of the chromatic polynomial $P_G(q)$
lie in the disc $|q| < C(\Delta)$,
where $C(\Delta)$ are semi-explicit constants (given by a variational formula)
satisfying $C(\Delta) \le 7.963907 \Delta$.\footnote{
   More recently, Borgs \cite{Borgs_06} has provided a
   simpler variational characterization
   of the constant $K = \lim_{\Delta\to\infty} C(\Delta)/\Delta
     \approx 7.963906$ than the one given by
   Sokal \cite[Proposition~5.4]{Sokal_01} ---
   compare eqs.~\reff{eq.sokalK} and \reff{eq.borgsK} below ---
   and Fern\'andez and Procacci \cite{FP_paper2} have provided,
   in an analogous way,
   a simpler variational characterization of the constants $C(\Delta)$.
   Furthermore, Fern\'andez and Procacci \cite{FP_paper2}
   have improved the constants $C(\Delta)$
   to smaller constants $C^*(\Delta)$,
   for which $K^* = \lim_{\Delta\to\infty} C^*(\Delta)/\Delta
     \approx 6.907652$.
} More generally, Sokal proved a bound on the zeros of the
multivariate Tutte polynomial (also known in statistical mechanics
as the Potts-model partition function, see
 \cite{Sokal_bcc2005,Potts_52,Wu_82,Wu_84})
\be
   Z_G(q, \w)   \;=\;
   \sum_{A \subseteq E}  q^{k(A)}  \prod_{e \in A}  w_e
 \label{def.ZG}
\end{equation}
[here $k(A)$ denotes the number of connected components
 in the subgraph $(V,A)$]
when the edge weights $\w = \{w_e\}$ lie in the
``complex antiferromagnetic regime'' $|1+w_e| \le 1$:

\begin{theorem}   {$\!\!\!$ \bf \protect\cite[Corollary~5.5]{Sokal_01} \ }
  \label{thm1.1}
Let $G=(V,E)$ be a loopless graph
equipped with complex edge weights $\w = \{ w_e \}_{e \in E}$
satisfying $|1 + w_e| \le 1$ for all $e$.
Then all the zeros of $Z_G(q,\w)$ lie in the disc
$|q| < K \Delta(G,\w)$,
where
\be
   \Delta(G,\w)   \;=\;  \max\limits_{x \in V} \sum\limits_{e\ni x}
   |w_e|
\ee
and
\begin{subeqnarray}
   K  & = &
   \min\left\{ L  \colon\;
                  \inf\limits_{\alpha > 0}
   \alpha^{-1} \sum\limits_{n=2}^\infty  e^{\alpha n} \, L^{-(n-1)} \,
       \frac{n^{n-1}}{n!}
   \;\le\;   1
       \right\}    \slabel{eq.sokalK} \\[2mm]
   & = &
   \min\limits_{a > 0} \,  {a+e^a \over \log(1+a e^{-a})}
        \slabel{eq.borgsK}  \\[2mm]
   & \approx & 7.963\:906\:075\:890\:002\:502\:\ldots \;.
\end{subeqnarray}
Moreover, we rigorously have $K \le 7.963907$.
\end{theorem}

\noindent Here the simpler formula \reff{eq.borgsK} for the constant
$K$ is due to Borgs \cite[Theorem~2.1]{Borgs_06}.

The purpose of this paper is to extend Sokal's bound
by removing the condition that $|1 + w_e| \le 1$ for all $e$.
More precisely, we shall prove:\footnote{
   A simpler but weaker version of this result
   can be found in the first and second preprint versions of this paper
   (\url{http://arxiv.org/abs/0810.4703v1} and \url{v2}).
}

\begin{theorem}
  \label{thm1.2c}
Let $G=(V,E)$ be a loopless graph equipped with complex edge weights
$\w = \{ w_e \}_{e \in E}$. Then all the zeros of $Z_G(q,\w)$ lie in
the disc
\be
   |q| \;<\; \widehat\scrk(\Psi(G,\w)) \: \widehat\Delta(G,\w)  \;,
 \label{e1.2c.0}
\ee
where
\begin{eqnarray}
   \widehat\Delta(G,\w)   & = &
      \max\limits_{x \in V} \sum\limits_{\substack{e\ni x \\ e=xy}}
          \min\left\{|w_{e}|,\, {|w_{e}|\over |1+w_{e}|}\right\}
 \prod\limits_{f \ni y}  \max\{1, |1+w_{f}|\}^{1/2} \label{e1.2c.1}
        \\[2mm]
   \Psi(G,\w)  & = &
      \max\limits_{x \in V} \prod\limits_{e\ni x}  \max\{1, |1+w_{e}|\}
   \quad \label{e1.2c.2}
\end{eqnarray}
and
\begin{subeqnarray}
   \widehat\scrk(\psi)  & = &
   \min\left\{ L  \colon\;
                  \inf\limits_{\alpha > 0} \,
   (e^\alpha-1)^{-1} \sum\limits_{n=2}^\infty  e^{\alpha n} \, \psi^{1/2} \,
         L^{-(n-1)} \, {n^{n-1} \over n!}
   \;\le\;   1
       \right\} \qquad     \slabel{def_kstar.a}\\[2mm]
   & = &
   \min\limits_{1 < y < 1+\psi^{-1/2}} \,
        {\psi^{-1/2} \, y \over (1+\psi^{-1/2}-y) \log y}
        \slabel{eq.FPKprimeb} \\[2mm]
   & = &
    \psi^{-1/2} \, W\biggl( \displaystyle {e \over 1+\psi^{-1/2}} \biggr)
     \bigg/
      \left[ 1 - W\biggl( \displaystyle {e \over 1+\psi^{-1/2}}
                  \biggr) \right]^2
        \slabel{eq.lambertKprimec} \\[2mm]
   & \le & 4\psi^{1/2} + 3   \;,
 \label{def_kstar}
 \label{e1.2c.6}
\end{subeqnarray}
where $W$ is the Lambert $W$ function \cite{Corless_96}, i.e.\ the
inverse function to $x \mapsto x e^x$.
\end{theorem}

When $|1 + w_e| \le 1$ for all $e$, we have $\widehat\Delta(G,\w) =
\Delta(G,\w)$ and $\Psi(G,\w) = 1$, so that Theorem~\ref{thm1.2c}
reduces in this case to Theorem~\ref{thm1.1} with an improved
constant \cite{FP_paper2} $K^* \equiv \widehat\scrk(1) =
W(e/2)/[1-W(e/2)]^2
 \approx 6.907\:651\:697\:774\:449\:218\:\ldots\;\,$.
This explicit formula
for the Fern\'andez--Procacci \cite{FP_paper2} constant $K^*$
appears to be new.

Let us also remark that the upper bound (\ref{def_kstar}d) gives
precisely the first two terms of the large-$\psi$ asymptotics of
$\widehat\scrk(\psi)$:  see equation \reff{eq.Kpsi.largepsi} in the Appendix.

Please note that both $\Psi(G,\w)$ and $\widehat\Delta(G,\w)$ involve
a {\em product}\/ over all edges incident to a given vertex rather
than a sum, and hence grow {\em exponentially}\/ (rather than
linearly) with the vertex degree whenever $|1+w_e| > 1$. The
resulting exponential dependence of the bound on $|q|$ given in
Theorem~\ref{thm1.2c} is not merely an artifact of our proof, but is
a genuine feature of the regime $|1+w_e| > 1$.\footnote{
   See also \cite[Remark~2 after Corollary~5.5]{Sokal_01}.
} To see this, it suffices to note that whenever one replaces an
edge~$e$ by $k$ edges in parallel, the effective couplings
$w_{e,{\rm eff}} = (1+w_e)^k - 1$ grow exponentially in $k$ when
$|1+w_e| > 1$ but only linearly when $|1+w_e| \le 1$. For instance,
the graph $G = K_2^{(k)}$ (a pair of vertices connected by $k$
parallel edges) with all edge weights equal has $Z_G(q,w) = q [q +
(1+w)^k - 1]$, so that we must take $|q| > |(1+w)^k - 1|$ to avoid a
root. This has roughly (but not exactly) the same dependence in $w$
and $k$ as the bound of Theorem~\ref{thm1.2c}. See
Example~\ref{exam.K2k} below for details.


When all edge weights are equal, the two factors
$\widehat\scrk(\Psi(G,\w))$ and $\widehat\Delta(G,\w)$
combine to produce a bound that grows linearly
with $\Psi(G,\w)$ as $\Psi(G,\w) \to\infty$. If we restrict
attention to {\em simple}\/ graphs, then with a little more
combinatorial work we can obtain a bound that grows only like
$\Psi(G,\w)^{1/2}$:

\begin{theorem}\label{thm1.3b}
Let $G=(V,E)$ be a simple graph (i.e.\ no loops or multiple edges)
equipped with complex edge weights $\w = \{ w_e \}_{e \in E}$. Then
all the zeros of $Z_G(q,\w)$ lie in the disc
\be
   |q| \;<\;  K^*_\mu  \:\Delta^*(G,\w) \;,
 \label{eq.thm1.3b.boundb}
\ee
where
\be
\Delta^*(G,\w)    =
      \max\limits_{x \in V} \sum\limits_{\substack{e \ni x \\ e=xy}}
          \min\left\{|w_{e}|,\, {|w_{e}|\over |1+w_{e}|^{1/2}}\right\}
 \prod\limits_{f \ni y}  \max\{1, |1+w_{f}|\}^{1/2}
      \label{def.Deltastarb}
\ee
and
$\mu = \widehat\Delta(G,\w)/\Delta^*(G,\w)$ and
\begin{subeqnarray}\label{e1.3b}
   K^*_\mu  & = &
  \min\left\{ L  \colon\,
                  \inf\limits_{\alpha > 0}
   (e^\alpha-1)^{-1} \sum\limits_{n=2}^\infty  e^{\alpha n}
        L^{-(n-1)} \, {[1+(n-1)\mu]^{n-2} \over (n-1)!}
   \le\;   1
       \right\} \qquad  \slabel{lemma.borgs.F0b.a}   \\[2mm]
   & = &
   \min\limits_{1 < y < 2} \,
        {y^\mu \over (2-y) \log y}
        \slabel{lemma.borgs.F0b.b} \\[3mm]
   & \le & 5  + 2\mu  \;.
      \slabel{lemma.borgs.F0b.c}
      \label{lemma.borgs.F0b}
\end{subeqnarray}
\end{theorem}

Please note that $0 < \mu \le 1$ because
$\min\{|w_{e}|,\, {|w_{e}|/|1+w_{e}|}\}
 \leq \min\{|w_{e}|,\, {|w_{e}|/|1+w_{e}|^{1/2}}\}$
for all $e\in E$, hence $\widehat\Delta(G,\w)\leq\Delta^*(G,\w)$.
The constant $K^*_\mu$ is an increasing function of $\mu \in (0,1]$,
but the variation is fairly weak:
we have $K^*_0 = W(2e) / [2 \, [W(2e)-1]^2]
   \approx 4.892888$
and $K^*_1 = K^* = W(e/2)/[1-W(e/2)]^2 \approx 6.907652$.
Thus, in the complex antiferromagnetic regime $|1+w_e| \le 1$ for all $e$,
where $\mu = 1$, Theorems~\ref{thm1.2c} and \ref{thm1.3b} give the same bound.


When $|1+w_e| > 1$, by contrast, Theorem~\ref{thm1.3b} is in most cases
a big improvement over Theorem~\ref{thm1.2c}:
this is because $K^*_\mu$ is always order 1
while $\widehat\scrk(\Psi(G,\w))$ is order $\Psi(G,\w)^{1/2}$.

Note that the bound \reff{e1.2c.0} involves a double maximum:
once over $x \in V$ in $\Psi(G,\w)$,
and once over $x \in V$ in $\widehat\Delta(G,\w)$.
Such a bound is ``unnatural'' in the sense that
if $G$ is a disjoint union $G = G_1 \uplus G_2$,
then the chromatic roots of $G$ are the union of those of $G_1$ and $G_2$,
and $\widehat\scrk(\Psi)$ and $\widehat\Delta$
are each the maximum of those for $G_1$ and $G_2$,
but the product $\widehat\scrk(\Psi) \, \widehat\Delta$ for $G$
{\em can exceed}\/ the maximum of those for $G_1$ and $G_2$
because one factor could be maximized for $G_1$ and the other for $G_2$
(see Example~\ref{exam.G1G2} below).
The bound \reff{eq.thm1.3b.boundb} has the virtue of avoiding
such a double maximum.
It is an open question whether a bound avoiding a double maximum
can be obtained for non-simple graphs.

On the other hand, in the bound \reff{eq.thm1.3b.boundb} we do pay a
price, compared to \reff{e1.2c.0}, by having
$\Delta^*(G,\w)$ in place of $\widehat\Delta(G,\w)$, since as noted above we have $\Delta^*(G,\w)\ge \widehat\Delta(G,\w)$.
In fact, the simple example $G=K_2$ shows that
the bound of Theorem~\ref{thm1.3b} can in some cases
be inferior to that of Theorem~\ref{thm1.2c},
by a factor of up to $K^*_0/4 \approx 1.223222$
(see Examples~\ref{exam.K2} and \ref{exam.thm1.2_beats_thm1.3} below).
But this seems to be the largest possible ratio of the two bounds.

It is curious that the bound of Theorem~\ref{thm1.3b}
is not always better than that of Theorem~\ref{thm1.2c},
despite using better ``ingredients'' in its proof;
the reasons for this will be discussed near the end of
Section~\ref{sec.final}.
It would be interesting to try to find a single natural bound that
simultaneously improves Theorems~\ref{thm1.2c} and \ref{thm1.3b}.

Please note also (see e.g.\ \cite{Sokal_bcc2005})
that if $G$ is a loopless graph with multiple edges,
then its multivariate Tutte polynomial
is identical to that of the underlying simple graph $\widehat G$
in which each set of parallel edges $e_1, \ldots, e_k$ in $G$
is replaced by a single edge $e$ in $\widehat G$ with weight
$\widehat w_{e} = \prod_{i=1}^k (1+w_{e_i}) - 1$.
So one is always free to apply Theorem~\ref{thm1.2c} or \ref{thm1.3b}
to $(\widehat G,\widehat \w)$
instead of applying Theorem~\ref{thm1.2c} to $(G,\w)$.
The following lemma concerning the behavior of
$\Psi(G,\w)$ and $\widehat\Delta(G,\w)$ under parallel reduction ---
which will be proven at the end of Section~\ref{sec.final} ---
implies that the bound we get by applying Theorem~\ref{thm1.2c}
to $(\widehat G,\widehat \w)$ will never be worse than
the bound we get by applying Theorem~\ref{thm1.2c} to $(G,\w)$.
So we can find our best bound for any given (multi)graph $G$
by constructing $(\widehat G,\widehat \w)$
and then taking the minimum of the bounds we obtain by applying
(\ref{e1.2c.0}) and (\ref{eq.thm1.3b.boundb})
to $(\widehat G,\widehat \w)$.

\begin{lemma}
   \label{lemma.parallel}
Let $w_1,w_2\in \mathbb C$ and put $w_3=(1+w_1)(1+w_2)-1$. Then
\begin{equation}
   \label{parallelprod}
\max\{1, |1+w_3|\}  \;\le\; \max\{1, |1+w_1|\} \, \max\{1, |1+w_2|\}
\end{equation}
and
\begin{equation}
   \label{parallelsum}
\min\left\{|w_3|,\, {|w_3|\over |1+w_3|}\right\}  \;\le\;
\min\left\{|w_1|,\, {|w_1|\over |1+w_1|}\right\} \,+\,
\min\left\{|w_2|,\, {|w_2|\over |1+w_2|}\right\}  \,.
\end{equation}
\end{lemma}


\bigskip

Sokal's proof of Theorem~\ref{thm1.1} involved the following steps:
\begin{itemize}
   \item[1.] Write the multivariate Tutte polynomial $Z_G(q,\w)$
      as the partition function of a polymer gas
      with weights depending on $q$ and $\w$
      (this is easy: see Section~\ref{sec.polymer} below).
   \item[2.] Invoke the Koteck\'y--Preiss \cite{Kotecky_86}
      condition for the nonvanishing of the partition function
      of a polymer gas.
   \item[3.] Control the polymer weights by bounding
      sums over connected subgraphs by sums over trees,
      using the Penrose inequality \cite{Penrose_67}.
      This step required $|1+w_e| \le 1$.
   \item[4.] Bound the total weight of $n$-vertex trees (or more generally,
      of connected subgraphs with $m$ edges) in $G$ that contain
      a specified vertex $x \in V$.
   \item[5.] Put everything together to prove that $Z_G(q,\w) \neq 0$
      whenever $q$ lies outside a specified disc.
\end{itemize}
Here we follow the same outline, but modify step 3
so as to allow arbitrary complex weights $w_e$.
In addition, in step 2 we replace the Koteck\'y--Preiss condition
by the more powerful Gruber--Kunz--Fern\'andez--Procacci
\cite{Gruber_71,FP} condition,
thereby slightly improving the numerical constant
along the lines of the work of Fern\'andez and Procacci \cite{FP_paper2}
for chromatic polynomials.
Finally, we need a slightly strengthened version
of the bound in step 4.

The plan of this paper is to treat each of these five steps
in successive sections.
Thus, in Section~\ref{sec.polymer} we recall how
the multivariate Tutte polynomial $Z_G(q,\w)$
can be written as the partition function of a polymer gas.
In Section~\ref{sec.convergence} we recall the
Koteck\'y--Preiss and Gruber--Kunz--Fern\'andez--Procacci conditions
for the nonvanishing of the partition function of a polymer gas.
In Section~\ref{sec.penrose} we recall the Penrose identity \cite{Penrose_67}
and show how to use it to bound the polymer weights
{\em without}\/ assuming that $|1+w_e| \le 1$;
this is our main new contribution.
In Section~\ref{sec.counting} we prove a bound
on the total weight of connected $m$-edge subgraphs in $G$
that contain a specified vertex~$x$;
this strengthens the bound of \cite{Sokal_01,maxmaxflow}
by taking specific account of the edges incident on $x$
and by introducing vertex weights.
In Section~\ref{sec.final} we put everything together
to prove Theorems~\ref{thm1.2c} and \ref{thm1.3b};
we also prove Lemma~\ref{lemma.parallel}.
Finally, in Section~\ref{sec.examples} we examine some examples
that shed light on the extent to which
Theorems~\ref{thm1.2c} and \ref{thm1.3b} are sharp or non-sharp.
In an Appendix we prove Lemma~\ref{lemma.oldcor.borgs} and some related facts.

\section{Polymer-gas representation of $\bm{Z_G(q,\w)}$}  \label{sec.polymer}
\zeq

In statistical mechanics, an {\em abstract polymer gas}\/
is a triple $(P, \xi, \mathcal{R})$
where $P$ is a finite
set (whose elements are called ``polymers"),
$\xi$ is a complex-valued function defined on $P$
(the value $\xi(p)$ is called the ``activity" or ``fugacity'' or ``weight''
 of the polymer $p \in P$),
and $\mathcal{R}\subseteq P\times P$  is a symmetric and reflexive relation
(called the ``incompatibility relation'').
Note that, since $\mathcal{R}$ is supposed reflexive,
we have $(p,p)\in \mathcal{R}$ for all $p\in P$.
Then the {\em partition function}\/ of the polymer gas $(P, \xi, \mathcal{R})$
--- a key quantity from which all thermodynamic properties of the system
can in principle be derived --- is defined by
\be
   \Xi(\xi)
   \;=\;
   \sum_{n=0}^\infty
   \sum\limits_{\substack{ \{p_1,\ldots,p_n\} \subseteq P  \\[0.7mm]
                           (p_i, p_j)\notin \mathcal{R} \: \forall i \neq j
                         }}
   \xi(p_1) \,\cdots\, \xi(p_n)
\ee
where the sum runs over unordered collections $\{p_1,\ldots,p_n\}$
of mutually compatible elements of $P$,
and the $n=0$ term in the sum is understood to contribute $1$.

In this section we recall how to rewrite the
multivariate Tutte polynomial $Z_G(q,\w)$ of a graph $G=(V,E)$
as the partition function of a polymer gas
living on the vertex set of $G$,
i.e.\ an abstract polymer gas whose polymers are nonempty subsets of $V$.
This easy result is due to Sokal and Kupiainen
\cite[Proposition 2.1]{Sokal_01}.

First, some notation:
If $H=(\sV,\sE)$ is a graph equipped with edge weights
$\w = \{w_e\}_{e \in \sE}$,
we denote by $C_H(\w)$ the generating polynomial
of connected spanning subgraphs of $H$, i.e.
\begin{equation}
   C_H(\w) \;=\;  \sum\limits_{\substack{
                                    A \subseteq \sE \\[0.7mm]
                                    (\sV,A) \, {\rm connected}
                               }}
                       \prod_{e \in A}  w_e
   \;.
 \label{def.CH}
\end{equation}
Note that $C_H(\w) \equiv 0$ if $H$ is disconnected.

If $G=(V,E)$ is a graph and $S \subseteq V$,
we denote by $G[S]$ the induced subgraph of $G$ on $S$,
i.e.\ $G[S]$ is the graph whose vertex set is $S$
and whose edges consist of all the edges of $G$
both of whose endpoints lie in $S$.

\begin{proposition}[polymer representation of the multivariate Tutte polynomial]
 \label{prop.polymer}
Let $G=(V,E)$ be a loopless graph
equipped with edge weights $\w = \{w_e\}_{e \in E}$.
Then
\be
   q^{-|V|} Z_G(q,\w)  \;=\;  \sum_{N=0}^\infty
      \;
      \sum_{\substack{
                \{S_1,\ldots,S_N\} \\[0.7mm]
                {\rm disjoint}
            }}
      \;
      \prod_{i=1}^n \xi(S_i)
      \;,
 \label{eq.polymer}
\ee
where the sum runs over unordered collections $\{S_1,\ldots,S_N\}$
of disjoint nonempty subsets of $V$,
and the weights $\xi(S)$ are given by
\be
   \xi(S)  \;=\;  \begin{cases}
                      q^{-(|S|-1)} C_{G[S]}(\w)  & \text{if } |S| \ge 2 \\[1mm]
                      0     & \text{if } |S| = 1
                  \end{cases}
 \label{def.polymer_weights}
\ee
[The $N=0$ term in the sum \reff{eq.polymer} is understood to contribute $1$.]
\end{proposition}

The identity \reff{eq.polymer} thus represents $q^{-|V|} Z_G(q,\w)$
as the partition function of a polymer gas
given by the triple $(P,\xi,\mathcal{R})$
with the polymer space $P$ being the set of all nonempty subsets of $V$,
the activity $\xi$  being the function defined in \reff{def.polymer_weights},
and the incompatibility relation $\mathcal{R}$ being nonempty intersection,
i.e.\ $(S,S')\in \mathcal{R}$ if and only if $S\cap S'\neq \emptyset$.
Note that, since the weight $\xi(S)$ vanishes for sets of cardinality 1
and also vanishes whenever the induced subgraph $G[S]$ is disconnected,
we can equivalently restrict our polymer set $P$
to be the set of all subsets $S \subseteq V$ of cardinality at least 2 and
for which $G[S]$ is connected.

Hereafter we will refer to a polymer gas in which polymers are subsets
of a given set $V$ and the incompatibility relation is nonempty intersection
as ``a gas of nonoverlapping polymers living on $V$".

\proofof{Proposition~\ref{prop.polymer}}
Starting from the definition \reff{def.ZG} of $Z_G(q,\w)$,
let us separate the terms in the sum according to the number $k$
of connected components [i.e.\ $k(A)=k$]
and according to the partition $\{S_1,\ldots,S_k\}$ of $V$
that is induced by the vertex sets of those connected components;
we will then sum over all ways of choosing edges within those
vertex sets $S_i$ so as to connect those vertices.
We thus have
\be
   Z_G(q,\w)  \;=\;  q^{|V|} \sum_{k \ge 1}
      \;
      \sum_{\substack{
                \{S_1,\ldots,S_k\}  \\[0.7mm]
                V = \biguplus S_i
            }}
      \;
      \prod_{i=1}^k q^{-(|S_i|-1)} C_{G[S_i]}(\w)
   \;,
\ee
where the sum runs over all unordered partitions $\{S_1,\ldots,S_k\}$ of $V$
into nonempty subsets, and we have used $|V| = \sum_{i=1}^k |S_i|$.
Note now that any set $S_i$ of cardinality~1
gets weight $q^{-(|S_i|-1)} C_{G[S_i]}(\w) = 1$
(here we have used the fact that $G$ is loopless).
So let us define $\{S'_1,\ldots,S'_N\}$
to be the subcollection of $\{S_1,\ldots,S_k\}$
consisting of the sets of cardinality $\ge 2$;
and let us note that there is a one-to-one correspondence between
unordered partitions $\{S_1,\ldots,S_k\}$ of $V$ into nonempty subsets
and unordered collections $\{S'_1,\ldots,S'_N\}$
of disjoint subsets of $V$ of cardinality at least 2
(which need not cover all of $V$:  indeed, the points not covered
 correspond to the singleton sets $S_i$ in the original partition).
Passing to $\{S'_1,\ldots,S'_N\}$ and dropping the primes,
we have \reff{eq.polymer}/\reff{def.polymer_weights}.
\qed

%
%
%
%
%
%
%
%
%

\section{Sufficient condition for the nonvanishing of a polymer-gas
    partition function}   \label{sec.convergence}
\zeq

Let $V$ be a finite set, and let $\{\rho(S)\}_{\emptyset \neq S \subseteq V}$
be a collection of complex weights associated to the nonempty subsets of $V$.
Consider now a gas of nonoverlapping polymers living on $V$,
with weights $\rho(S)$:  the partition function of such a polymer gas is,
by definition,
\be
   \Xi  \;=\;
      \sum_{N=0}^\infty \:
      \sum_{\substack{
                \{S_1,\ldots,S_N\} \\[0.7mm]
                {\rm disjoint}
            }}
      \;
      \prod_{i=1}^N \rho(S_i)
      \;,
 \label{eq.polymer.general}
\ee
where the sum runs over unordered collections $\{S_1,\ldots,S_N\}$
of disjoint nonempty subsets of $V$,
and the $N=0$ term in \reff{eq.polymer.general}
is understood to contribute $1$.
The following proposition ---
essentially proven almost four decades ago
by Gruber and Kunz \cite[Section~4, cf.~eq.~(33)]{Gruber_71}
but largely forgotten,
and then rediscovered very recently by Fern\'andez and Procacci
\cite[eq.~(3.17)]{FP} with a new proof ---
gives a sufficient condition for the nonvanishing
of a polymer-gas partition function:

\begin{proposition}[Gruber--Kunz--Fern\'andez--Procacci condition]
   \label{prop.GKFP}
Let $V$ be a finite set, and let $\{\rho(S)\}_{\emptyset \neq S \subseteq V}$
be complex weights associated to the nonempty subsets of $V$.
Suppose that there exists a number $\alpha > 0$ such that
\be
   \sup_{x \in V} \sum_{S \ni x}  e^{\alpha |S|} \, |\rho(S)|  \;\le\;
      e^\alpha - 1
   \;.
 \label{eq.GKFP}
\ee
Then
\be
   \Xi  \;\equiv\;
   \sum_{N=0}^\infty \:
   \sum_{\substack{
                \{S_1,\ldots,S_N\} \\[0.7mm]
                {\rm disjoint}
            }}
      \;
      \prod_{i=1}^n \rho(S_i)
   \;\,\neq\;\, 0  \;.
\ee
\end{proposition}

\noindent
See also \cite{Bissacot_10} for an extremely simple proof
of Proposition~\ref{prop.GKFP} by induction on $V$.

In the slightly less powerful Koteck\'y--Preiss \cite{Kotecky_86} condition,
the term $e^\alpha - 1$ on the right-hand side of \reff{eq.GKFP}
is replaced by $\alpha$.

\bigskip

{\bf Remark.}
Suppose that (as happens in all nontrivial cases)
there exists a set $S$ with $|S| \ge 2$ and $\rho(S) \neq 0$.
Then the hypothesis that there exists $\alpha > 0$
such that \reff{eq.GKFP} holds can be rewritten as
\be
   \inf_{\alpha > 0} (e^\alpha - 1)^{-1}
      \sup_{x \in V} \sum_{S \ni x}  e^{\alpha |S|} \, |\rho(S)|  \;\le\  1
   \;,
 \label{eq.GKFP.bis}
\ee
since in this case the infimum on the left-hand side of \reff{eq.GKFP.bis}
will always be attained at some $\alpha > 0$.\footnote{
   If there exists a set $S$ with $|S| \ge 2$ and $\rho(S) \neq 0$,
   then the function $f(\alpha)$
   being minimized on the left-hand side of \reff{eq.GKFP.bis}
   is a continuous function that tends to $+\infty$
   as $\alpha \downarrow 0$ and as $\alpha \uparrow \infty$,
   hence its minimum is attained.

   There is one exceptional case in which \reff{eq.GKFP.bis} holds
   but there does not exist $\alpha > 0$ such that \reff{eq.GKFP} holds:
   namely, if $\rho(S) = 0$ whenever $|S| \ge 2$
   and in addition we have $\max\limits_{x \in V} |\rho(\{x\})| = 1$.
   Indeed, if $\rho(S) = 0$ for $|S| \ge 2$,
   we have $\Xi = \prod\limits_{x \in V} [1 + \rho(\{x\})]$,
   which vanishes when at least one $\rho(\{x\})$ equals $-1$;
   so \reff{eq.GKFP.bis} {\em fails}\/ (barely) to imply
   $\Xi \neq 0$ in this case.
}
We will use the Gruber--Kunz--Fern\'andez--Procacci condition
in the form \reff{eq.GKFP.bis}.

%
%

\section{A bound on $\bm{C_H(\w)}$ via the Penrose identity}
  \label{sec.penrose}
\zeq

In this section we recall the Penrose identity \cite{Penrose_67} and
show how it can be used to bound a sum over connected subgraphs by a
sum over trees {\em even in the absence of the hypothesis $|1+w_e|
\le 1$}\/.

Let $H=(\sV,\sE)$ be a graph. Recall that $C_H(\w)$ denotes the
generating polynomial of connected spanning subgraphs of $H$:
\begin{equation}
   C_H(\w) \;=\;  \sum\limits_{\substack{
                                    A \subseteq \sE \\[0.7mm]
                                    (\sV,A) \, {\rm connected}
                               }}
                       \prod_{e \in A}  w_e
   \;.
\end{equation}
We denote by $T_H(\w)$ the generating polynomial of spanning trees
in $H$:
\begin{equation}
   T_H(\w) \;=\;  \sum\limits_{\substack{
                                    A \subseteq \sE \\[0.7mm]
                                    (\sV,A) \, {\rm tree}
                               }}
                       \prod_{e \in A}  w_e
   \;.\label{tree}
\end{equation}

Let $\scrc$ (resp.\ $\scrt$) be the set of subsets $A \subseteq \sE$
such that $(\sV,A)$ is connected (resp.\ is a tree). Clearly $\scrc$
is an increasing family of subsets of $\sE$ with respect to
set-theoretic inclusion, and the minimal elements of $\scrc$ are
precisely those of $\scrt$ (i.e.\ the spanning trees). It is a
nontrivial combinatorial fact
--- apparently first discovered by Penrose \cite{Penrose_67} ---
that the (anti-)complex $\scrc$ is {\em partitionable}\/: that is,
there exists a map ${\bf R} \colon\, \scrt \to \scrc$ such that
${\bf R}(T) \supseteq T$ for all $T \in \scrt$ and $\scrc =
\biguplus_{T \in \scrt} [T, \, {\bf R}(T)]$ (disjoint union), where
$[E_1,E_2]$ denotes the Boolean interval $\{ A \colon\; E_1
\subseteq A \subseteq E_2 \}$.
We call any such map~${\bf R}$ a~{\em partition scheme}\/.
In fact, many alternative choices of ${\bf R}$ are available\footnote{
   See for example\ \cite{Penrose_67},
   \cite[Sections 7.2 and 7.3]{Bjorner_92},
   \cite[Section 8.3]{Ziegler_95},
   \cite[Sections 2 and 6]{Gessel_96},
   \cite[Proposition 13.7 et seq.]{Biggs_93},
   \cite[Proposition 4.1]{Sokal_01}
   and \cite[Lemma~2.2]{Scott-Sokal}.
},
and most of our arguments will not depend on any specific choice
of ${\bf R}$. An immediate consequence of the existence of ${\bf R}$
is the following simple but fundamental identity:

\begin{proposition}[Penrose identity \protect\cite{Penrose_67}]
Let ${\bf R} \colon\, \scrt \to \scrc$
be any partition scheme.
Then
\begin{subeqnarray}
      C_H(\w)
      & = &
\sum\limits_{\substack{
                  T \subseteq \sE  \\[0.7mm]
                  (\sV,T) \, {\rm tree}
             }}
      \; \prod\limits_{e \in T} w_e
      \sum\limits_{T \subseteq A \subseteq {\bf R}(T)}
      \; \prod\limits_{e \in A \setminus T} w_e
         \\[3mm]
      & = &
\sum\limits_{\substack{
                  T \subseteq \sE  \\[0.7mm]
                  (\sV,T) \, {\rm tree}
             }}
      \; \prod\limits_{e \in T} w_e
      \; \prod\limits_{e \in {\bf R}(T) \setminus T} (1 + w_e)
      \;.
    \label{penrose_identity}
    \slabel{penrose_identity_b}
\end{subeqnarray}
\end{proposition}

If $|1+w_e| \le 1$ for all $e$, then it is obvious that we can take
absolute values everywhere in \reff{penrose_identity_b} and drop the
factors $|1+w_e|$, yielding:

\begin{proposition}[Penrose inequality \protect\cite{Penrose_67}]
     \label{prop.penrose}
Let $H=(\sV,\sE)$ be a graph equipped with complex edge weights
$\w = \{ w_e \}_{e \in \sE}$ satisfying $|1 + w_e| \le 1$ for all $e$.
Then
\be
      |C_H(\w)|  \;\le\;  T_H(|\w|)   \;.
  \label{eq.prop.penrose}
\ee
\end{proposition}

{\bf Remark.}  By using a specific choice of the map ${\bf R}$
(namely, that of Penrose \cite{Penrose_67}), Fern\'andez and
Procacci \cite{FP} have recently shown how to improve
Proposition~\ref{prop.penrose} when $w_e \in \{-1,0\}$ for all $e$;
and this improvement plays a key role in their proof of the
Gruber--Kunz--Fern\'andez--Procacci condition
(Proposition~\ref{prop.GKFP}) for polymer gases with hard-core
repulsive interactions. See also Fern\'andez {\em et al.}\/
\cite{FKSU} for a generalization to $-1 \le w_e \le 0$, which leads
to an improved convergence criterion for the Mayer expansion in
lattice gases with soft repulsive interactions. \qed

Let us now show what can be done {\em without}\/ the hypothesis
$|1 + w_e| \le 1$.
Given a vertex $x$ in a graph $H = (\sV,\sE)$,
we denote by $\sE(x)$
the set of edges of $H$ incident on $x$.
For any subset $A \subseteq \sE$, let us write
\begin{subeqnarray}
   A_+  & = &   \{ e \in A \colon\:  |1+w_e| > 1 \}   \\[1mm]
   A_-  & = &   \{ e \in A \colon\:  |1+w_e| \le 1 \}
\end{subeqnarray}

\begin{proposition}[extended Penrose inequality]
     \label{prop.extpenrose}
Let $H=(\sV,\sE)$ be a loopless
graph equipped with complex edge weights $\w = \{ w_e \}_{e \in \sE}$.
Then
\begin{subeqnarray}
   |C_H(\w)|
   & \le &
   T_H(|\w'|) \, \prod_{e \in \sE} \max\{1, |1+w_{e}|\}
      \slabel{eq.extpenrose.a}
       \\[1mm]
   & = &
   T_H(|\w'|) \, \prod_{y \in \sV} \prod_{e \in \sE(y)}
                     \max\{1, |1+w_{e}|\}^{1/2}
      \slabel{eq.extpenrose.b}
      \label{eq.extpenrose}
\end{subeqnarray}
where
\begin{equation}
   w'_e
   \;=\;
   \begin{cases}
        w_e                & \text{\rm if } |1+w_e| \le 1 \\[3mm]
        \displaystyle {w_e \over 1+w_e}  & \text{\rm if } |1+w_e| > 1
   \end{cases}
 \label{def.weprime}
\end{equation}
\end{proposition}

Note that if $|1 + w_e| \le 1$ for all $e$,
then $\w' = \w$ and $\max\{1, |1+w_{e}|\} = 1$ for all~$e$,
so Proposition~\ref{prop.extpenrose} is a genuine extension
of Proposition~\ref{prop.penrose}.

\proofof{Proposition~\ref{prop.extpenrose}}
In the Penrose identity \reff{penrose_identity_b},
multiply and divide the summand by
\hbox{$\prod\limits_{e \in T_+} (1+w_e)$:}
this yields
\be
   C_H(\w)
   \;=\;
   \sum\limits_{\substack{
                  T \subseteq \sE  \\[0.7mm]
                  (\sV,T) \, {\rm tree}
             }}
      \; \prod\limits_{e \in T} w'_e
      \; \prod\limits_{e \in ({\bf R}(T) \setminus T) \,\cup\, T_+} (1 + w_e)
      \;.
\ee
Taking absolute values and using the trivial bound
\be
   \prod\limits_{e \in ({\bf R}(T) \setminus T) \,\cup\, T_+} |1 + w_e|
   \;\le\;
   \prod_{e \in \sE} \max\{1, |1+w_{e}|\}
   \;,
 \label{eq.proof.extpenrose}
\ee
we obtain \reff{eq.extpenrose.a}.
Then \reff{eq.extpenrose.b} follows by observing that
each edge $e \in \sE$ is incident on precisely two vertices
(since $H$ is loopless).
\qed

{\bf Remark.}
Quite a lot has been thrown away in \reff{eq.proof.extpenrose}.
Can we do better in a usable way?
\qed

\bigskip

If we assume that the graph $H$ is {\em simple}\/ (i.e.\ has no
loops or multiple edges), then we can get a slightly better bound:

\begin{proposition}[extended Penrose inequality for simple graphs]
     \label{prop.extpenrose.simple}
Let $H=(\sV,\sE)$ be a {\em simple}\/ graph
(i.e.\ no loops or multiple edges)
equipped with complex edge weights $\w = \{ w_e \}_{e \in \sE}$.
Then, for any vertex $x \in \sV$, we have
\begin{subeqnarray}
   |C_H(\w)|
   & \le &
   T_H(|\bwdoubleprimex|) \,
   \prod_{e \in \sE \setminus \sE(x)} \max\{1, |1+w_{e}|\}
      \slabel{eq.extpenrose.simple.a}
       \\[1mm]
   & \le &
   T_H(|\bwtildex|) \,
      \prod_{y\in V\setminus \{x\}} \,
      \prod_{e\in \sE(y)}\max\{1, |1+w_{e}|\}^{1/2}
      \slabel{eq.extpenrose.simple.b}
      \label{eq.extpenrose.simple}
\end{subeqnarray}
where
\begin{equation}
   \wdoubleprimex_e
   \;=\;
   \begin{cases}
        w_e                & \text{\rm if } |1+w_e| \le 1
                                \text{\rm\ or } e \in \sE(x)   \\[3mm]
        \displaystyle {w_e \over 1+w_e}
                           & \text{\rm if } |1+w_e| > 1 \text{\rm\ and }
                             e \in \sE \setminus \sE(x)   \\[3mm]
   \end{cases}
 \label{def.wedoubleprime}
\end{equation}
and
\begin{equation}
   \wtildex_e
   \;=\;
   \begin{cases}
        w_e                & \text{\rm if } |1+w_e| \le 1 \\[3mm]
        \displaystyle {w_e \over |1+w_e|^{1/2}}
                           & \text{\rm if } |1+w_e| > 1 \text{\rm\ and }
                             e \in \sE(x)   \\[6mm]
        \displaystyle {w_e \over |1+w_e|}
                           & \text{\rm if } |1+w_e| > 1 \text{\rm\ and }
                             e \in \sE \setminus \sE(x)   \\[3mm]
   \end{cases}
 \label{def.wetilde}
\end{equation}
\end{proposition}




Please note that \reff{eq.extpenrose.simple.b}
is indeed an improvement of \reff{eq.extpenrose.b},
because the product $\prod_{e \in \sE(x)} \max\{1, |1+w_{e}|\}^{1/2}$
more than compensates the factors
$|\wtildex_e/w'_e| = \max\{1, |1+w_{e}|\}^{1/2}$
for the {\em subset}\/ of edges in $\sE(x)$ that happen to lie in
any given spanning tree~$T$.

The proof of Proposition~\ref{prop.extpenrose.simple}
will be based on the following key combinatorial fact (to be proven later):

\begin{lemma}
   \label{lemma.specialR}
Let $H=(\sV,\sE)$ be a {\em simple}\/ graph and let $x \in \sV$ be any vertex.
Then there exists a partition scheme ${\bf R}$ with the property that
${\bf R}(T) \setminus T$ does not contain any edge incident on $x$.
\end{lemma}

%

\bigskip\noindent
{\sc Proof of Proposition~\ref{prop.extpenrose.simple},
  assuming Lemma~\ref{lemma.specialR}.\ }
In the Penrose identity \reff{penrose_identity_b},
multiply and divide the summand by
\hbox{$\prod\limits_{e \in [T \setminus \sE(x)]_+} (1+w_e)$:}
this yields
\be
   C_H(\w)
   \;=\;
   \sum\limits_{\substack{
                  T \subseteq \sE  \\[0.7mm]
                  (\sV,T) \, {\rm tree}
             }}
      \; \prod\limits_{e \in T} \wdoubleprimex_e
      \; \prod\limits_{e \in  [{\bf R}(T) \setminus T] \,\cup\,
                              [T \setminus \sE(x)]_+} (1 + w_e)
      \;.
 \label{eq.proof.extpenrose.simple.1}
\ee
Choosing the partition scheme as in Lemma~\ref{lemma.specialR},
we have ${\bf R}(T) \setminus T \subseteq \sE \setminus \sE(x)$
and hence
\be
   \prod\limits_{e \in [{\bf R}(T) \setminus T] \,\cup\, [T \setminus \sE(x)]_+}
      |1 + w_e|
   \;\le\;
   \prod_{e \in \sE \setminus \sE(x)} \max\{1, |1+w_{e}|\}
   \;.
 \label{eq.proof.extpenrose.simple.2}
\ee
Taking absolute values in \reff{eq.proof.extpenrose.simple.1}
and using \reff{eq.proof.extpenrose.simple.2}, we obtain
\be
   |C_H(\w)|
   \;\le\;
   \sum\limits_{\substack{
                  T \subseteq \sE  \\[0.7mm]
                  (\sV,T) \, {\rm tree}
             }}
      \; \prod\limits_{e \in T} |\wdoubleprimex_e|
      \; \prod_{e \in \sE \setminus \sE(x)} \max\{1, |1+w_{e}|\}
      \;,
 \label{eq.proof.extpenrose.simple.3}
\ee
which is \reff{eq.extpenrose.simple.a}.

Now observe that
\be
   \prod_{e \in \sE \setminus \sE(x)} \max\{1, |1+w_{e}|\}
   \;=\;
   { \prod\limits_{y\in V \setminus x} \, \prod\limits_{e\in \sE(y)}
         \max\{1, |1+w_{e}|\}^{1/2}
     \over
     \prod\limits_{e \in \sE(x)} \max\{1, |1+w_{e}|\}^{1/2}
   }
 \label{eq.proof.extpenrose.simple.4}
\ee
since the numerator of \reff{eq.proof.extpenrose.simple.4} counts
every edge in $\sE \setminus \sE(x)$ twice and every edge in $\sE(x)$ once.
If in the denominator of \reff{eq.proof.extpenrose.simple.4}
we replace the product over $e \in \sE(x)$ by the smaller product
over $e \in \sE(x) \cap T$, we get an upper bound;
inserting this into \reff{eq.proof.extpenrose.simple.3}
yields \reff{eq.extpenrose.simple.b}.
\qed

Let us conclude this section by proving Lemma~\ref{lemma.specialR}.
This proof --- unlike all the preceding results in this section ---
depends on a specific choice of the map ${\bf R}$, namely the one
used by Penrose in his original paper \cite{Penrose_67}.
Let us briefly recall Penrose's construction (see \cite{FP,FKSU} for more
details). We assume that $H=(\sV,\sE)$ is a {\em simple}\/ graph,
and we choose (arbitrarily) an ordering of the vertex set $\sV$ by
numbering the vertices $1,2,\ldots,n$ (where $n=|\sV|$). We consider
the vertex 1 to be the root, and denote it by $r$. If $T \subseteq
\sE$ is the edge set of a spanning tree in $H$ [that is, $(\sV,T)$
is a tree], then for each $x \in \sV$ we denote by ${\rm dist}_T(x)$
the graph-theoretic distance in the tree $(\sV,T)$ from the root~$r$
to the vertex~$x$. Given $T$, the vertex set $\sV$ is thus
partitioned into ``generations'', defined as the sets of vertices at
a given distance from the root~$r$.

The Penrose map ${\bf R} \colon\, T \mapsto {\bf R}(T)$ is then
defined as follows. For any tree $T \subseteq \sE$, the edge set
${\bf R}(T) \supseteq T$ is obtained from $T$ by adjoining all edges
$e \in \sE$ that either
\begin{itemize}
   \item[(a)]  connect two vertices in the same generation
      [i.e.\ at equal distance from the root~$r$ in the tree $(\sV,T)$ ---
       note that no such edge can belong to $T$], or
   \item[(b)]  connect a vertex $x$ to a vertex $x'$ in the preceding
      generation [i.e.\ with ${\rm dist}_T(x') = {\rm dist}_T(x) - 1$]
      that is higher-numbered than the parent of $x$
      [here the parent of $x$ is the unique vertex $y$
       with ${\rm dist}_T(y) = {\rm dist}_T(x) - 1$ such that $xy \in T$].
\end{itemize}
It can be shown \cite{Penrose_67,FP,FKSU} that ${\bf R}$ is indeed a
partitioning map in the sense that $\scrc$ is the disjoint union of
Boolean intervals $[T, \, {\bf R}(T)]$. Furthermore, it follows
immediately from this construction that ${\bf R}(T) \setminus T$
cannot contain any edge incident on the root~$r$; that is, ${\bf
R}(T) \setminus T \subseteq \sE \setminus \sE(r)$.\footnote{
  We remark that this would no longer be the case in a generalization
  to the Penrose construction to non-simple graphs.
  In such a generalization, we would also order the edges
  connecting each pair of vertices, and we would add to the definition of
  ${\bf R}(T)$ a third case:
  \begin{itemize}
      \item[(c)] connect a vertex $x$ to its parent $y$
          by any edge that is higher-numbered
          than the edge connecting $x$ to $y$ in $T$.
  \end{itemize}
  We would then no longer be able to guarantee that
  ${\bf R}(T) \setminus T$ contains no edges incident on the root~$r$;
  rather, we could assert only that
  ${\bf R}(T) \setminus T$ cannot contain any edge incident on the root~$r$
  that is the {\em lowest-numbered}\/ among its set of parallel edges.
}
Since any vertex could have been chosen as the root,
Lemma~\ref{lemma.specialR} is proven.

\bigskip

{\bf Remark.}  Lemma~\ref{lemma.specialR} suggests the following
combinatorial question:
Let $H=(\sV,\sE)$ be a graph (simple or not).
For which subsets $S \subseteq \sE$ does there exist a partition scheme
${\bf R}$ with the property that
${\bf R}(T) \setminus T \subseteq \sE \setminus S$ for all $T$?
The same question can also be posed for matroids.
\qed

%
%
%
%
%
%
%

\section{Bounds on connected $\bm{m}$-edge subgraphs containing a
      specified vertex}
   \label{sec.counting}
\zeq

In this section consider a loopless graph $G=(V,E)$
equipped with nonnegative real edge weights $\{w_e\}_{e \in E}$
{\em and}\/ nonnegative real vertex weights $\{w_v\}_{v \in V}$.
Let us define the weighted sum over connected subgraphs
$G' = (V',E') \subseteq G$
that contain a specified vertex $x$ and have exactly $m$ edges:
\be
   c_m(x;G,\w)  \;=\;  \sum_{\substack{
                            G' = (V',E') \subseteq G \\[0.7mm]
                            G' \ {\rm connected} \\[0.7mm]
                            V' \ni x \\[0.7mm]
                            |E'| = m
                        }}
                  \;
                  \prod_{e \in E'} w_e \prod_{v \in V'} w_v
   \;,
\ee
where we write $\w = \{ w_e \}_{e \in E}\cup \{ w_v \}_{v \in V}$.
We will abbreviate $c_m(x;G,\w)$ to $c_m(x)$ when it is obvious
which weighted graph $(G,\w)$ we are referring to.
Now define the weighted degree at $x$ by
\be
   d(x;G,\w) \;=\; \sum_{e=xy\in E} w_e \, w_y
\ee
(note that this contains a factor $w_y$ for each edge $e=xy$ incident to $x$ but {\em not}\/ a factor $w_x$),
and define the maximum weighted degree by
\be
   \Delta(G,\w)   \;=\;  \max\limits_{x \in V}  d(x;G,\w)
   \;.
\ee

The following bound on $c_m(x)$ extends an earlier result of the
third author \cite[Proposition 4.5]{Sokal_01}, which is obtained by
putting $w_v=1$ for all $v\in V$ and using the fact that both
$d(x;G,\w)$ and $\Delta(G-x,\w|_{G-x})$ are bounded above by
$\Delta(G,\w)$.

\begin{proposition}
   \label{prop.counting.aldo}
Let $G=(V,E)$ be a loopless graph equipped with nonnegative real weights
$\w = \{ w_e \}_{e \in E}\cup \{ w_v \}_{v \in V}$, and let $x\in V$.
Suppose that either $w_v \geq 1$ for all $v\in V$ or $G$ is simple.
Then
\be
   c_m(x)
   \;\le\;
   { w_x \: d(x;G,\w) \: [d(x;G,\w)+m\Delta(G-x,\w|_{G-x})]^{m-1}  \over m!}
     \;
 \label{eq.counting.aldo0}
\ee
for all $m\geq 0$.
\end{proposition}

We remark that the bound \reff{eq.counting.aldo0}
need not hold if we remove the hypothesis that
either $w_v \geq 1$ for all $v\in V$ or $G$ is simple.
Consider, for instance, the graph $G = K_2^{(m)}$
consisting of two vertices $x,y$ joined by $m \ge 2$ parallel edges.
Put $w_x=w_y=w$ and $w_e=1$ for all $e\in E$.
Then $c_m(x)=w^2$,
while the right-hand side of \reff{eq.counting.aldo0}
is $m^m w^{m+1}/m!$, which is less than $c_m(x)$ when $w$ is small enough.

\bigskip

In the proof of Proposition~\ref{prop.counting.aldo} it will be
convenient to employ the quantities \be
   C(m,\kappa)   \;=\;
   \begin{cases}
      \kappa (m+\kappa)^{m-1} / m!     & \text{for } m \ge 1 \\
      1                                & \text{for } m = 0
   \end{cases}
\ee defined for integer $m \ge 0$ and real $\kappa$. Then
\reff{eq.counting.aldo0} can be rewritten in the form \be
   c_m(x)   \;\le\;  w_x \: C(m,d/\Delta) \: \Delta^m   \;
    \label{eq.counting.aldo1}
\ee where $d = d(x;G,\w)$ and $\Delta = \Delta(G-x, \w|_{\,G-x})$.

Our proof of Proposition~\ref{prop.counting.aldo} uses induction on
$m$, and
is similar to the first proof of  \cite[Proposition
7.1]{maxmaxflow}. It relies on the following properties of
$C(m,\kappa)$:
\begin{enumerate}
\item  For each integer $m \ge 0$, $C(m,\kappa)$ is a polynomial of degree $m$
in $\kappa$, with nonnegative coefficients.  In particular,
$C(m,\kappa)$ is an increasing function of $\kappa$ for real $\kappa
\ge 0$.

\item  Generating function:  If ${\sf C}(z)$ solves the equation
\be
   {\sf C}(z)   \;=\;  e^{z {\sf C}(z)}
   \;,
 \label{gen_fn_eqn}
\ee then \be
   {\sf C}(z)^\kappa \;=\;  \sum\limits_{m=0}^\infty C(m,\kappa) \, z^m
 \label{gen_fn}
\ee for all real $\kappa$; this follows from the Lagrange inversion
formula. Moreover, the series \reff{gen_fn} is absolutely convergent
for $|z| \le 1/e$ and satisfies ${\sf C}(1/e) = e$.

\item For integer $k \ge 1$,
\be
   C(m,k)  \;=\;  \!\!\!\sum
   _{\substack{
                  m_1, \ldots, m_k \ge 0 \\[0.7mm]
                  m_1 + \cdots + m_k = m }}
          \prod\limits_{i=1}^k C(m_i,1)
   \;.
     \label{cmk.identity1}
\ee This is an immediate consequence of \reff{gen_fn}.

\item For all real $\kappa$ and $z$,
\be
   C(m,\kappa)  \;=\;  \sum\limits_{f=0}^m {z^f \over f!} \, C(m-f,\kappa-z+f)
   \;.
     \label{cmk.identity2}
\ee See \cite[eq.~(7.7)]{maxmaxflow}.
\end{enumerate}

For any subset $F \subseteq E$, we use the notation
$w(F) = \prod_{e \in F} w_e$.
Also, for any $F \subseteq E(x)$, we denote by $Y^F$
the set of vertices of $V-x$ that are incident with edges in $F$,
and we write $j(F) = |Y^F|$ for the number of such vertices.
Please observe that $j(F)\le |F|$; and if the graph $G$ is simple,
then $j(F)= |F|$.

Our proof of Proposition~\ref{prop.counting.aldo}
will be based on the following two lemmas:

\begin{lemma}\label{recursion}
Let $G=(V,E)$ be a loopless graph equipped with nonnegative real weights
$\w = \{ w_e \}_{e \in E}\cup \{ w_v \}_{v \in V}$, and let $x \in V$.
For each $F\subseteq E(x)$, let $Y^F =
\{x_1^F,x_2^F,\ldots,x_{j(F)}^F\}$ be a labeling of the vertices of
$V-x$ that are incident with edges in $F$. Then, for all $m\geq 1$,
\be c_m(x;G,\w) \le w_x\sum\limits_{\emptyset\neq F \subseteq E(x)}
   w(F) \, \sum_{\substack{
                  m_1, \ldots, m_{j(F)} \ge 0 \\[0.7mm]
                  m_1 + \cdots + m_{j(F)} = m-|F|}}
           \! \prod\limits_{i=1}^{j(F)} c_{m_i}(x_i^F;G-x,\w|_{G-x}) \;.
\ee
  \label{eq.lemma.recursion}
\end{lemma}

\proof
Similar to that given for Facts~1 and 2 in \cite[Section~7]{maxmaxflow}.
\qed

%

\begin{lemma}  {$\!\!\!$ \bf \protect\cite[Lemma~7.2]{maxmaxflow} \ }
  \label{subsetweight}
Let $S$ be a set in which each element $e\in S$ is given a
nonnegative real weight $w_e$.
Then, for each integer $f \ge 0$, we have \be
   \sum_{\substack{
                               F \subseteq S \\[0.7mm]
                                       |F| = f }}
w(F)
   \;\le\;
   {1 \over f!} \, \left( \sum\limits_{e \in S} w_e \right) ^{\! f}
   \;.
\ee
\end{lemma}

\proofof{Proposition~\ref{prop.counting.aldo}} Let $d=d(x;G,\w)$ and
$\Delta=\Delta(G-x,\w|_{G-x})$. We will prove
\reff{eq.counting.aldo0}/\reff{eq.counting.aldo1} by induction on
$m$. The statement holds trivially when $m=0$, so let us assume that
$m\geq 1$. By Lemma \ref{recursion},
\begin{eqnarray}
   c_m(x)
   & \le &
   w_x \sum\limits_{\emptyset\neq F \subseteq E(x)}
   w(F) \, \sum_{\substack{
                  m_1, \ldots, m_{j(F)} \ge 0 \\[0.7mm]
                  m_1 + \cdots + m_{j(F)} = m-|F|}}
           \! \prod\limits_{i=1}^{j(F)} c_{m_i}(x_i^F;G-x,\w|_{G-x})
      \nonumber \\[2mm]
   & \leq &
   w_x \sum\limits_{\emptyset\neq F \subseteq E(x)}
   w(F) \, \sum_{\substack{
                  m_1, \ldots, m_{j(F)} \ge 0 \\[0.7mm]
                  m_1 + \cdots + m_{j(F)} = m-|F|}}
           \! \prod\limits_{i=1}^{j(F)} w_{x_i^F} \; C(m_i, \, 1) \; \Delta^{m_i}
    \nonumber \\[2mm]
   & = &
   w_x \sum\limits_{\emptyset\neq F \subseteq E(x)}
    C(m-|F|, \, j(F)) \; \Delta^{m-|F|}\; w(F) \;\prod\limits_{i=1}^{j(F)}  w_{x_i^F}
    \nonumber \\[2mm]
   & \le &
   w_x \;\sum\limits_{f=1}^m C(m-f, f) \;
\Delta^{m-f}\;\sum_{\substack{F \subseteq E(x)\\[0.7mm]
|F|=f}}\;
   \prod\limits_{e=xx_i^F\in F} w_e \: w_{x_i^F}
\end{eqnarray}
where the second line used the induction hypothesis
\reff{eq.counting.aldo0} applied to the graph $G-x$ (note that $m_i
< m$)
and the fact that $d(v;G-x,\w|_{G-x})\leq \Delta$ for all $v\in V-x$;
the third line used the identity \reff{cmk.identity1};
and the last line used $j(F)\leq |F|$, the fact that $C(m,k)$ is an
increasing function of $k$, and the hypothesis that either
$w_{x_i^F} \geq 1$ for all $1\leq i\leq j(F)$ or $G$ is simple.
Using Lemma \ref{subsetweight}, we have
\begin{eqnarray}
   c_m(x)  & \le &
      w_x \,\Delta^m\,\sum\limits_{f=1}^m  {(d/\Delta)^f \over f!} \, C(m-f,f)
    \nonumber \\[2mm]
   & = &
      w_x \,\Delta^m\,\sum\limits_{f=0}^m  {(d/\Delta)^f \over f!} \, C(m-f,f)
    \nonumber \\[2mm]
   & = &  w_x \,C(m,d/\Delta) \, \Delta^m   \;,
\end{eqnarray}
where the second line used $C(m,0) = 0$ for $m \ge 1$, and the last
line used identity \reff{cmk.identity2} with $\kappa=z=d/\Delta$.
This proves (\ref{eq.counting.aldo1}). 
\qed

We now combine Proposition~\ref{prop.counting.aldo}
with the extended Penrose inequalities from Section~\ref{sec.penrose}:

\begin{proposition}
   \label{prop.aldo.new}
Let $G=(V,E)$ be a loopless graph equipped with complex edge weights
$\w = \{ w_e \}_{e \in E}$.
Let $x\in V$ and let $n$ be a positive integer. Then

\be
   \sum_{\substack{
                  S \ni x \\[0.7mm]
                  S\subseteq V\\[0.7mm]
                  |S|=n}}
   |C_{G[S]}(\w)|
   \;\le\;
  {n^{n-1} \over n!} \, \widehat\Delta(G,\w)^{n-1}\;
    \prod_{e\in E(x)}\max\{1,|1+w_e|\}^{1/2}
 \label{eq.aldo.new.0}
\ee
where $\widehat\Delta(G,\w)$ is defined in \reff{e1.2c.1}.
Furthermore, if $G$ is simple, then
\be
   \sum_{\substack{
                  S \ni x \\[0.7mm]
                  S\subseteq V\\[0.7mm]
                  |S|=n}}
  |C_{G[S]}(\w)|
  \;\le\;
  {\Delta^*(G,\w)\over (n-1)!}
\left[ \Delta^*(G,\w)+(n-1)\widehat\Delta(G,\w) \right]^{n-2}
 \label{eq.aldo.new.1}
\ee
where $\Delta^*(G,\w)$ is defined in \reff{def.Deltastarb}.
\end{proposition}

\proof
We first prove (\ref{eq.aldo.new.0}). Construct a nonnegative real weight
function $\widehat\w$ on $V\cup E$ by putting $\widehat w_y=\prod_{e\in
E(y)}\max\{1,|1+w_e|\}^{1/2}$ for all $y\in V$, and $\widehat w_e= |w'_e|$
for all $e\in E$,
where $w'_e$ is defined in \reff{def.weprime}.
For $y\in S\subseteq V$ let $E(y;G[S])$ denote the set of edges of
$G[S]$ incident on $y$.
By bound \reff{eq.extpenrose.b} of Proposition~\ref{prop.extpenrose}, we have
\begin{subeqnarray}
\sum_{\substack{
                S \ni x \\[0.7mm]
                S\subseteq V\\[0.7mm]
                |S|=n}}
   |C_{G[S]}(\w)|
   &\leq&
 \sum_{\substack{
                S \ni x \\[0.7mm]
                S\subseteq V\\[0.7mm]
                |S|=n}}
   T_{G[S]}(|\w'|) \,
  \prod_{y\in S} \: \prod_{e\in E(y;G[S])}\max\{1,|1+w_e|\}^{1/2} \slabel{eq.extprop5.a} \quad \\
  &\leq&
  c_{n-1}(x;G,\widehat\w)
\end{subeqnarray}
since the $n$-vertex trees are a subset of the connected graphs with
$n-1$ edges, and $E(y;G[S]) \subseteq E(y)$.
Inequality  (\ref{eq.aldo.new.0}) now follows by applying
Proposition~\ref{prop.counting.aldo}, using the fact that
$d(x;G,\widehat \w)$ and $\Delta(G-x,\widehat \w|_{G-x})$ are both
bounded above by $\Delta(G,\widehat \w)=\widehat\Delta(G,\w)$.

We next prove (\ref{eq.aldo.new.1}). Construct a weight function
$\w^*$ on $V\cup E$ by putting $w_x^*=1$, $w_y^*=\prod_{e\in
E(y)}\max\{1,|1+w_e|\}^{1/2}$ for all $y\in V\setminus \{x\}$, and
$w_e^*= |\wtildex_e|$ for all $e\in E$,
where $\wtildex_e$ is defined in \reff{def.wetilde}.
By bound \reff{eq.extpenrose.simple.b} of
Proposition~\ref{prop.extpenrose.simple}, we have
\begin{subeqnarray}
\sum_{\substack{
                S \ni x \\[0.7mm]
                S\subseteq V\\[0.7mm]
                |S|=n}}
|C_{G[S]}(\w)|
&\leq&
 \sum_{\substack{
                S \ni x \\[0.7mm]
                S\subseteq V\\[0.7mm]
                |S|=n}}
T_{G[S]}(|\bwtildex|)\,\prod_{y\in S\setminus\{x\}} \, \prod_{e\in E(y;G[S])}
   \max\{1,|1+w_e|\}^{1/2} \slabel{eq.extprop5.b}  \qquad \\
&\leq& c_{n-1}(x;G,\w^*)
\end{subeqnarray}
by the same reasoning as before.
Inequality  (\ref{eq.aldo.new.1}) now follows by applying
Proposition~\ref{prop.counting.aldo}, using the facts that
$d(x;G,\w^*)\leq {\Delta^*(G,\w)}$ and
$\Delta(G-x,\w^*|_{G-x})\leq \widehat\Delta(G,\w)$.
\qed

\section{Proof of Theorems~\ref{thm1.2c} and \ref{thm1.3b}
    and Lemma~\ref{lemma.parallel}}   \label{sec.final}
\zeq

We can now put together the results of the preceding sections to
prove Theorems~\ref{thm1.2c} and \ref{thm1.3b}. At the end of this
section we will also prove Lemma~\ref{lemma.parallel}.

We begin by stating an analytic lemma that will be needed in proving
the equivalence between the various versions
(\ref{e1.2c.6}a--d) and (\ref{lemma.borgs.F0b}a--c)
of our bounds.
To avoid disrupting the flow of the argument,
the proof of this lemma is deferred to an Appendix.

\begin{lemma}
   \label{lemma.oldcor.borgs}
For $\lambda \ge 0$ and $\beta > 0$, define the function
\be
   F_\lambda(\beta)   \;=\;
   \min\left\{ L  \colon\;
                  \inf\limits_{\alpha > 0} \,
   (e^\alpha-1)^{-1} \sum\limits_{n=2}^\infty  e^{\alpha n} \,
         L^{-(n-1)} \, {[1 + (n-1)\lambda]^{n-2} \over (n-1)!}
   \;\le\;   \beta
       \right\} \;. \qquad
\ee
Then
\be
   F_\lambda(\beta)
   \;=\;
   \min\limits_{1 < y < 1+\beta} \,
        {\beta y^{\lambda} \over (1+\beta-y) \log y}
   \;.
\ee
Moreover,
\be
   F_1(\beta)
   \;=\;
     \beta\,W\biggl( \displaystyle {e \over 1+\beta} \biggr)
     \bigg/
      \left[ 1 - W\biggl( \displaystyle {e \over 1+\beta}
                  \biggr) \right]^2
\ee
where $W$ is the Lambert $W$ function \cite{Corless_96},
 i.e.\ the inverse function to $x \mapsto x e^x$.
Finally, for $0 \le \lambda \le 1$ we have
\be
   F_\lambda(\beta)  \;\le\;  4\beta^{-1} + (1+2\lambda)
   \;.
\ee
\end{lemma}

\proofof{Theorem~\ref{thm1.2c}}
We want to show that $Z_G(q,\w) \neq 0$ whenever
$|q| \ge \widehat{\scrk}(\Psi(G,\w)) \, \widehat{\Delta}(G,\w)$.
We will do this by verifying the condition \reff{eq.GKFP.bis}
for the polymer weights \reff{def.polymer_weights}, which we recall are
\be
   \xi(S) \;=\; q^{-(|S|-1)} C_{G[S]}(\w)  \qquad\hbox{for } |S| \ge 2
   \;.
 \label{eq.sec.final.polymer_weights}
\ee
By inequality \reff{eq.aldo.new.0} of Proposition~\ref{prop.aldo.new},
for each $x\in V$ and each $n \ge 1$ we have
\begin{subeqnarray}
   \sum_{\substack{
                  S \ni x \\[0.7mm]
                  S\subseteq V\\[0.7mm]
                  |S|=n}}|C_{G[S]}(\w)|   \;&\le&\;{n^{n-1} \over n!} \,
                  \widehat\Delta(G,\w)^{n-1}\;\prod_{e\in
                  E(x)}\max\{1,|1+w_e|\}^{1/2}   \\[-3mm]
                  &\le\;&{n^{n-1} \over n!} \,
                  \widehat\Delta(G,\w)^{n-1}\;\Psi(G,\w)^{1/2}  \;.
 \label{eq.aldo.new.3}
\end{subeqnarray}
Therefore, the condition \reff{eq.GKFP.bis} for the weights
\reff{def.polymer_weights}/\reff{eq.sec.final.polymer_weights} is
verified as soon as
\begin{eqnarray}
   \inf_{\alpha > 0} (e^\alpha - 1)^{-1}
   \sum\limits_{n=2}^\infty  e^{\alpha n} \, |q|^{-(n-1)} \,
         {n^{n-1} \over n!} \, \widehat\Delta(G,\w)^{n-1} \, \Psi(G,\w)^{1/2}  &  \le
  1  \;.
 \label{eq.finalproof.finalcondition}
\end{eqnarray}
If we set $L= |q| \, \widehat\Delta(G,\w)^{-1}$ and $\psi = \Psi(G,\w)$
in \reff{eq.finalproof.finalcondition},
this is precisely the inequality contained in the
right-hand side of \reff{def_kstar.a}.
So $Z_G(q,\w) \neq 0$ whenever $L \ge \widehat\scrk({\Psi(G,\w)})$,
i.e.\ whenever $|q| \ge\widehat\scrk({\Psi(G,\w)}) \,\widehat\Delta(G,\w)$,
where $\widehat\scrk(\psi)$ is defined by \reff{def_kstar.a}.
The equivalence of \reff{def_kstar.a} with (\ref{def_kstar}b,c)
and the inequality (\ref{def_kstar}d)
follow from Lemma~\ref{lemma.oldcor.borgs}
once we observe that $\widehat\scrk(\psi) = F_1(\psi^{-1/2})$.
\qed

\medskip

\proofof{Theorem~\ref{thm1.3b}} We modify the proof of
Theorem~\ref{thm1.2c} by using (\ref{eq.aldo.new.1}) in place of
(\ref{eq.aldo.new.0}).

Since $G$ is simple, it follows from \reff{eq.aldo.new.1} that
for each $x\in V$ and each $n \ge 1$ we have
\begin{eqnarray}
   \sum_{\substack{
                  S \ni x \\[0.7mm]
                  S\subseteq V\\[0.7mm]
                  |S|=n}}|C_{G[S]}(\w)|
   \;&\le&\; {\Delta^*(G,\w)\over (n-1)!}
\left[ \Delta^*(G,\w)+(n-1)\widehat\Delta(G,\w) \right]^{n-2}
\nonumber \\[-4mm]
& = &\; \Delta^*(G,\w)^{n-1}  \;
           {[1+(n-1)\mu]^{n-2}  \over (n-1)!}
 \label{eq.thm1.3proof.sumTGSb}
\end{eqnarray}
where $\mu = \widehat\Delta(G,\w)/\Delta^*(G,\w)$. Therefore,
the condition \reff{eq.GKFP.bis} for the weights
\reff{def.polymer_weights}/\reff{eq.sec.final.polymer_weights} is
verified as soon as
\be
   \inf_{\alpha>0}  (e^\alpha-1)^{-1}
      \sum_{n\ge 2} e^{\alpha n} \,
        [|q|^{-1}  \Delta^*(G,\w)]^{n-1}  \:
        {[1+(n-1)\mu]^{n-2} \over (n-1)!}
   \;\le\; 1  \;.
 \label{eq.proof.thm1.3bb}
\ee
If we set $L= |q| \, \Delta^*(G,\w)^{-1}$,
this is precisely the inequality contained
in the right-hand side of \reff{lemma.borgs.F0b.a}.
So $Z_G(q,\w) \neq 0$ whenever $L \ge K^*_\mu$,
i.e.\ whenever $|q| \ge K^*_\mu \,  \Delta^*(G,\w)$,
where $K^*_\mu = F_\mu(1)$ is defined by \reff{lemma.borgs.F0b.a}.
The equivalence of \reff{lemma.borgs.F0b.a} with \reff{lemma.borgs.F0b.b}
and the inequality \reff{lemma.borgs.F0b.c}
then follow from Lemma~\ref{lemma.oldcor.borgs}.
\qed

\bigskip



{\bf Discussion.}
1.  We can now understand why the apparently minor improvement from
\reff{eq.extpenrose.b} to \reff{eq.extpenrose.simple.b}
leads to the significant improvement (in most cases) of the final bound
from Theorem~\ref{thm1.2c} to Theorem~\ref{thm1.3b},
namely, replacing a growth $\sim \Psi(G,\w)^{1/2}$ by $1$.
Indeed, we can see using Lemma~\ref{lemma.oldcor.borgs} that whenever we have a bound
of the form
\be
   \sum_{\substack{
             S \ni x \\[0.7mm]
             |S| = n
         }}
   |C_{G[S]}(\w)|
   \;\le\;  {[1+\lambda(n-1)]^{n-2} \over (n-1)!} \, D^{n-1} \, \Psi^{b}
   \;,
\ee
we will obtain a bound on the roots of $Z_G(q,\w)$ of the form
\be
   |q|  \;<\;  D \, F_\lambda(\Psi^{-b})
  \;.
\ee
The bound \reff{eq.extpenrose.b} gives rise to inequality \reff{eq.aldo.new.0},
which in turn allows us to deduce Theorem~\ref{thm1.2c}
by taking $D=\widehat\Delta$, $\lambda=1$ and $b = 1/2$.
On the other hand, the bound \reff{eq.extpenrose.simple.b}
gives inequality \reff{eq.aldo.new.1}, which allows us to deduce
Theorem~\ref{thm1.3b} by taking $D=\Delta^*$, $\lambda=\widehat\Delta/\Delta^*$ and $b = 0$.

\medskip

2.  Let us compare the bounds provided by Theorems~\ref{thm1.2c}
and \ref{thm1.3b}:
\begin{subeqnarray}
   \hbox{Theorem~\ref{thm1.2c}:}
   & &
   \widehat\scrk(\Psi(G,\w)) \: \widehat\Delta(G,\w)  \\[1mm]
   \hbox{Theorem~\ref{thm1.3b}:}
   & &
   K^*_\mu \: \Delta^*(G,\w)
\end{subeqnarray}
where $\mu = \widehat\Delta(G,\w)/\Delta^*(G,\w) \in (0,1]$.
Their ratio is therefore
\be
   {\hbox{Theorem~\ref{thm1.3b}} \over \hbox{Theorem~\ref{thm1.2c}}}
   \;=\;
   {K^*_\mu\:  \over \mu \: \widehat\scrk(\Psi(G,\w))}
   \;=\;
   {F_\mu(1)  \over \mu \: F_1(\Psi(G,\w)^{-1/2})}
   \;.
\ee
Now, it is not difficult to see that
${\Delta}^*(G,\w) \leq \widehat\Delta(G,\w) \, \Psi(G,\w)^{1/2}$,
or in other words $\Psi(G,\w)^{-1/2} \le \mu$.\footnote{
   {\sc Proof.}   For each edge $e=xy$ we have
   \begin{eqnarray*}
      \min\left\{|w_{e}|,\, {|w_{e}|\over |1+w_{e}|^{1/2}}\right\}
      & = &
      \min\left\{|w_{e}|,\, {|w_{e}|\over |1+w_{e}|}\right\}
         \:\times\:  \max\{1, |1+w_{e}|\}^{1/2}   \\
      & \le &
      \min\left\{|w_{e}|,\, {|w_{e}|\over |1+w_{e}|}\right\}
         \:\times\:  \Psi(G,\w)^{1/2}
      \;.
   \end{eqnarray*}
   Multiplying this by $\prod_{f \ni y}  \max\{1, |1+w_{f}|\}^{1/2}$, summing  over $e \ni x$, and taking the maximum over $x \in V$,
   we obtain the desired inequality.
}
Since $F_1(\beta)$ is a decreasing function of $\beta$
(see Proposition~\ref{prop.borgs}(a) in the Appendix),
we have $F_1(\Psi(G,\w)^{-1/2}) \ge F_1(\mu)$ and hence
\be
   {\hbox{Theorem~\ref{thm1.3b}} \over \hbox{Theorem~\ref{thm1.2c}}}
   \;\le\;
   {F_\mu(1)  \over  \mu \, F_1(\mu)}
   \;\equiv\; g(\mu)
   \;.
\ee
Both $F_\mu(1)$ and $\mu \, F_1(\mu)$
are increasing functions of $\mu$ [see Proposition~\ref{prop.borgs}(a,b)],
but their ratio $g(\mu)$ does not have any obvious monotonicity.
{\em Numerically}\/ we find that $g(\mu)$
decreases from the value $K^*_0/4 \approx 1.223222$ at $\mu=0$
to a minimum value $\approx 0.930714$ at $\mu \approx 3.70249$,
and then increases to 1 as $\mu\to\infty$.
We have not succeeded in {\em proving}\/ that $g(\mu) \le g(0)$
for $\mu \in [0,1]$, but if is true we can conclude that
Theorem~\ref{thm1.3b} is never more than a factor
$\approx 1.223222$ worse than Theorem~\ref{thm1.2c}.
In any case we have
\be
   g(\mu)  \;\le\;
   {F_1(1) \over \lim\limits_{\mu \to 0} \mu \, F_1(\mu)}
   \;=\;
   {K^*_1 \over 4}
   \;\approx\;
   1.726913
   \quad\hbox{for } \mu \in [0,1]
   \;.
\ee
We shall see in Examples~\ref{exam.K2} and \ref{exam.thm1.2_beats_thm1.3}
that Theorem~\ref{thm1.3b} can indeed be up to a factor
$\approx 1.223222$ worse than Theorem~\ref{thm1.2c}.

\medskip

3.  It is curious that the bound of Theorem~\ref{thm1.3b}
is not always better than that of Theorem~\ref{thm1.2c},
despite using a better ``ingredient'' in its proof:
namely, the bound \reff{eq.extpenrose.simple.b}
from Proposition~\ref{prop.extpenrose.simple}
always beats the bound~\reff{eq.extpenrose.b}
from Proposition~\ref{prop.extpenrose}.
How is it that the final result can sometimes be worse?

The explanation is that the ratio of the bounds
\reff{eq.extpenrose.simple.b} and \reff{eq.extpenrose.b}
\be
    {\hbox{\reff{eq.extpenrose.simple.b}}
     \over
     \hbox{\reff{eq.extpenrose.b}}
    }
    \;=\;
    { T_H(|\bwtildex|)
      \over
      T_H(|\w'|) \,\prod_{e\in \sE(x)}\max\{1, |1+w_{e}|\}^{1/2}
    }
\ee
is the product of a a ``good'' factor $\prod_{e\in \sE(x)}\max\{1, 
|1+w_{e}|\}^{-1/2}$ and a ``bad'' factor $ T_H(|\bwtildex|)/T_H(\w')$.
Now, the ``bad'' factor $T_H(|\bwtildex|)/T_H(|\w'|)$
is always bounded by $\prod_{e\in \sE(x)}\max\{1, |1+w_{e}|\}^{1/2}$
--- which is why \reff{eq.extpenrose.simple.b} is always better than
  \reff{eq.extpenrose.b} ---
so it follows that
\be
    {    \sum\limits_{S \ni x, \, |S| = n}  T_{G[S]}(|\bwtildex|)
      \over
         \sum\limits_{S \ni x, \, |S| = n}  T_{G[S]}(|\w'|)
    }
    \;\le\;
    \prod_{e\in \sE(x)}\max\{1, |1+w_{e}|\}^{1/2}\le  \Psi(G,\w)^{1/2}
    \;.
  \label{eq.ratio.sumT}
\ee
But there is no guarantee that the {\em upper bounds}\/
on the numerator and denominator of \reff{eq.ratio.sumT},
obtained by applying respectively the bounds
\reff{eq.aldo.new.3} and \reff{eq.thm1.3proof.sumTGSb},
will also have a ratio $\le \Psi(G,\w)^{1/2}$.
Indeed, it can happen that this {\em fails}\/
(see Examples~\ref{exam.K2} and \ref{exam.thm1.2_beats_thm1.3}).

It is, nevertheless, somewhat disconcerting that
Theorem~\ref{thm1.3b} is not always better than Theorem~\ref{thm1.2c}.
It would be nice to find a single natural bound that
simultaneously improves Theorems~\ref{thm1.2c} and \ref{thm1.3b}.
\qed

\bigskip

Finally, let us prove Lemma~\ref{lemma.parallel}
concerning the behavior of
$\Psi(G,\w)$ and $\widehat\Delta(G,\w)$ under parallel reduction:

\proofof{Lemma~\ref{lemma.parallel}}
Inequality (\ref{parallelprod}) follows immediately from the fact that
$(1+w_1)(1+w_2)=1+w_3$.
To prove (\ref{parallelsum}), let us consider the following cases:
\medskip

{\em Case 1}\/: $|1+w_1|\leq 1$ and $|1+w_2|\leq 1$. Then
$\min\left\{|w_i|,\, {|w_i|\over |1+w_i|}\right\}=|w_i|$ for
$1\leq i\leq 3$, so we just have to prove that $|w_3| \le |w_1|+|w_2|$.
Since $w_3 = w_1+w_2+w_1w_2$, we have
\begin{eqnarray}
 && |w_3| \;=\; |w_1+w_2+w_1w_2| \;=\; |w_1+w_2(1+w_1)|
    \;\le\; |w_1|+|w_2(1+w_1)|
    \qquad\qquad
    \nonumber \\
   & & \hspace*{2.5in} \;=\;  |w_1|+|w_2|\,|1+w_1| \;\le\; |w_1|+|w_2|
    \qquad\qquad
\end{eqnarray}
since $|1+w_1|\leq 1$.

\medskip

{\em Case 2}\/: $|1+w_1|\geq 1$ and $|1+w_2|\geq 1$. Then
$\min\left\{|w_i|,\, {|w_i|\over |1+w_i|}\right\}={|w_i|\over |1+w_i|}$ for
$1\leq i\leq 3$. Let $w_i'=-{w_i\over 1+w_i}$ for $1\leq i\leq 3$,
so that $1+w'_i = (1+w_i)^{-1}$ for $1\leq i\leq 3$
and hence $(1+w'_1)(1+w'_2)=1+w'_3$.
Since $|1+w'_1|\leq 1$ and $|1+w'_2|\leq 1$,
we may apply Case 1 to $w'_1, w'_2, w'_3$ to deduce that
$|w'_3| \le |w'_1|+|w'_2|$, as required.

\medskip

{\em Case 3}\/: $|1+w_1|\leq 1$, $|1+w_2|\geq 1$ and $|1+w_1|\,|1+w_2|\leq 1$.
Then
$\min\left\{|w_i|,\, {|w_i|\over |1+w_i|}\right\}=|w_i|$ for $i\in\{1,3\}$,
and $\min\left\{|w_2|,\, {|w_2|\over |1+w_2|}\right\}={|w_2|\over |1+w_2|}$.
By hypothesis we have $|1+w_1|\leq|1+w_2|^{-1}$. Hence
\be
   |w_3| \;=\; |w_1+w_2(1+w_1)| \;\le\;  |w_1| + |w_2|\,|1+w_1|
     \;\le\; |w_1|+{|w_2|\over |1+w_2|} \;,
\ee
as required.

\medskip

{\em Case 4}\/: $|1+w_1|\leq 1$, $|1+w_2|\geq 1$ and $|1+w_1|\,|1+w_2|\geq 1$.
Then
$\min\left\{|w_1|,\, {|w_1|\over |1+w_1|}\right\}=|w_1|$,  and
$\min\left\{|w_i|,\, {|w_i|\over |1+w_i|}\right\}={|w_i|\over |1+w_i|}$
for $i\in\{2,3\}$.
Let $w_i'=-{w_i\over 1+w_i}$ for $1\leq i\leq 3$.
Then $|1+w'_1|\geq 1$ and $|1+w'_2|\leq 1$
with $|1+w'_1|\,|1+w'_2|\leq 1$,
so we may apply Case 3 (with indices 1 and 2 interchanged) to deduce that
$|w'_3| \le {|w'_1|\over |1+w'_1|}+|w'_2| = |w_1| +|w'_2|$,
as required.
\qed

{\bf Remark.}
We suspect that the transformation
\be
   w'  \;=\;  - \, {w \over 1+w}
 \label{eq.FAFduality}
\ee
employed in Cases 2 and 4,
which satisfies $(1+w') = (1+w)^{-1}$
and hence preserves the parallel-connection law
$(1+w_1)(1+w_2) = 1+w_3$, may have other applications
in the study of the multivariate Tutte polynomial.
This transformation is involutive [i.e.\ $(w')'=w$],
maps the complex antiferromagnetic regime $|1+w|\le 1$
onto the complex ferromagnetic regime $|1+w'|\ge 1$
and vice versa,
and maps the real antiferromagnetic regime $-1 \le w \le 0$
onto the real ferromagnetic regime $0 \le w' \le +\infty$
and vice versa.
In the physicists' notation $w = e^J - 1$
where $J$ is the Potts-model coupling,
the transformation \reff{eq.FAFduality} takes the simple form $J' = -J$,
which makes its properties obvious.

\section{Examples}  \label{sec.examples}
\zeq

In this section we examine some examples that shed light on
the extent to which Theorems~\ref{thm1.2c} and \ref{thm1.3b}
are sharp or non-sharp.
For each weighted graph $(G,\w)$, we attempt to compute or estimate the quantity
\be
   Q_{\rm max}(G,\w)  \;=\;
   \max\{ |q| \colon\; Z_G(q,\w) = 0 \}
\ee
and compare it to the upper bounds given by Theorem~\ref{thm1.2c}
and Theorem~\ref{thm1.3b}.
In what follows we abbreviate
$\widehat\Delta(G,\w)$, $\Delta^*(G,\w)$, $\Psi(G,\w)$, $Q_{\rm max}(G,\w)$
by $\widehat\Delta$, $\Delta^*$, $\Psi$, $Q_{\rm max}$.

\begin{example}
   \label{exam.K2}
\rm
Let $G = K_2$, where the single edge has weight $w$.
Then $Z_{K_2}(q,w) = q(q+w)$, so that $Q_{\rm max} = |w|$.
On the other hand, if $|1+w| \ge 1$ we have
$\widehat\Delta = |w|/|1+w|^{1/2}$, ${\Delta}^* = |w|$,
$\Psi = |1+w|$ and $\mu = \widehat\Delta/{\Delta}^* = 1/|1+w|^{1/2}$.
Theorem~\ref{thm1.2c} gives the bound
$|q| < \widehat\scrk(\Psi) \, \widehat\Delta$,
which behaves like $4|w|$ as $|w| \to \infty$,
while Theorem~\ref{thm1.3b} gives the bound
$|q| < K^*_\mu \, {\Delta}^*$,
which behaves like $K_0^* |w| \approx 4.892888 |w|$
as $|w| \to \infty$.
So Theorem~\ref{thm1.2c} is off by a factor of 4 from the truth,
while Theorem~\ref{thm1.3b} is off by a factor of $\approx 4.892888$
from the truth.
In particular, Theorem~\ref{thm1.3b} is worse than Theorem~\ref{thm1.2c}
by a factor tending to $K_0^*/4 \approx 1.223222$.

For the special case of $G=K_2$,
the convergence conditions \reff{eq.finalproof.finalcondition}
and \reff{eq.proof.thm1.3bb},
which were used in the proofs of Theorems~\ref{thm1.2c}
and \ref{thm1.3b}, respectively, become
\begin{eqnarray}
   \inf_{\alpha > 0} (e^\alpha - 1)^{-1}
   e^{2\alpha} \, |q|^{-1} \, \widehat\Delta(G,\w) \, \Psi(G,\w)^{1/2}
   & \le &
   1
       \label{eq.finalproof.finalconditionb}  \\[2mm]
   \inf_{\alpha>0}  (e^\alpha-1)^{-1}
      e^{2\alpha} \, |q|^{-1} \, \Delta^*(G,\w) \,
   & \le &
   1
       \label{eq.proof.thm1.3b?}
\end{eqnarray}
because the only polymer in the graph $K_2$ has size $n=2$.
Since $\widehat\Delta(G,\w) \, \Psi(G,\w)^{1/2} =
       \Delta^*(G,\w) = |w|$,
we have
\be
   \hbox{\reff{eq.finalproof.finalconditionb}}
   \;\Longleftrightarrow\;
   \hbox{\reff{eq.proof.thm1.3b?}}
   \;\Longleftrightarrow\;
   |q| \,\ge\, 4|w|
   \;,
 \label{eq.polymers.K2}
\ee
which differs from the truth $Q_{\rm max} = |w|$ by a factor of 4.
We can understand this behavior as follows:

1) The lost factor of 4 comes from the fact that,
for a polymer gas consisting of a single polymer $S$
of cardinality $|S|=2$,
the Gruber--Kunz--Fern\'andez--Procacci condition
(Proposition~\ref{prop.GKFP})
gives $\Xi \neq 0$ whenever $|\rho(S)| \le 1/4$,
whereas the truth is that $\Xi \neq 0$ whenever $|\rho(S)| < 1$.

2) Though the convergence condition \reff{eq.finalproof.finalcondition}
involves a sum $\sum_{n=2}^\infty$,
the terms for $n > 2$ make a negligible contribution
in the limit $|w| \to \infty$ because
\be
   |q|^{-(n-1)} \, \widehat\Delta(G,\w)^{n-1} \, \Psi(G,\w)^{1/2}
   \;=\;
   (|w|/|q|)^{n-1} \, |1+w|^{-(n-2)/2}
   \;,
\ee
which tends to zero as $|w| \to \infty$ whenever
$|q| \ge {\rm const} \times |w|$ and $n > 2$.
That is why Theorem~\ref{thm1.2c} is off from the truth
by the {\em same}\/ factor 4 that we see in \reff{eq.polymers.K2},
despite the fact that its proof allows for arbitrarily large polymers
that do not occur when $G=K_2$.

3) By contrast, in the convergence condition \reff{eq.proof.thm1.3bb},
the terms with $n > 2$ do {\em not}\/ disappear
in the limit $|w| \to \infty$ with $|q|$ of order $|w|$, because
\be
   [|q|^{-1} \,  \Delta^*(G,\w) ]^{n-1}
   \;=\;
   (|w|/|q|)^{n-1}
\ee
is of order 1 for all $n$.
This is why Theorem~\ref{thm1.3b} is off from the truth
by {\em more}\/ than the factor 4 that we see in \reff{eq.polymers.K2};
we lose an additional factor $K_0^*/4 \approx 1.223222$
by allowing for nonexistent large polymers.
\qed
\end{example}

\begin{example}
   \label{exam.thm1.2_beats_thm1.3}
\rm
In {\em any}\/ simple graph $G$ with at least one edge,
we can choose weights $\w$ such that Theorem~\ref{thm1.2c}
beats Theorem~\ref{thm1.3b} by a factor
arbitrarily close to $K_0^*/4 \approx 1.223222$.
It suffices to take $w_e = w$ (with $|1+w| \ge 1$)
on all the edges of a nonempty matching,
and $w_e = w_0$ on all other edges;
then as $w_0 \to 0$ we have
$Q_{\rm max} \to |w|$,
$\widehat\Delta \to |w|/|1+w|^{1/2}$, ${\Delta}^* \to |w|$,
$\Psi \to |1+w|$ and $\mu= \widehat\Delta/{\Delta}^* \to 1/|1+w|^{1/2}$.
So the comparison of the bounds is the same as for $G=K_2$,
and Theorem~\ref{thm1.2c} beats Theorem~\ref{thm1.3b} by a factor
tending to $K_0^*/4 \approx 1.223222$ as $|w| \to\infty$.

For instance,
let $G$ be the $n$-cycle $C_n$ with $n \ge 3$,
taking $w_e = w$ for exactly one edge and $w_e=w_0$ for all other edges.
Then $Z_G(q,w) = (q+w)(q+w_0)^{n-1} + w w_0^{n-1} (q-1)$.
As $|w| \to \infty$ at fixed $n$ and $w_0$,
we have $Q_{\rm max}(G,\w) = |w| + o(|w|)$.
On the other hand, if $|1+w_0| \ge 1$ and $|w| \gg |w_0|$
we have $\widehat\Delta(G,\w) = |w_0|+|1+w_0|^{1/2}|w|/|1+w|^{1/2}$,
$\Delta^*(G,\w) =  |w_0||1+w_0|^{1/2} + |1+w_0|^{1/2}|w|$
and $\Psi(G,\w) = |1+w_0| \, |1+w|$.
Therefore, as $|w| \to \infty$ the bounds of
Theorems~\ref{thm1.2c} and \ref{thm1.3b}
are
$4 |1+w_0| |w| + O(|w|^{1/2})$
and $K^*_0 |1+w_0|^{1/2}|w| + O(1)$, respectively,
where $K^*_0  \approx 4.892888$.
Both of these bounds have the correct order of magnitude
as $|w| \to \infty$ at fixed $n$ and $w_0$,
but are off by a constant factor
($4 |1+w_0|$ or $K^*_0 |1+w_0|^{1/2}$, respectively).
The bound given by Theorem~\ref{thm1.2c} is better
than that given by Theorem~\ref{thm1.3b} when $|1+w_0|$ is small,
and worse when $|1+w_0|$ is large.
\qed
\end{example}

\begin{example}
   \label{exam.K2k}
\rm
Let $G = K_2^{(k)}$ (a pair of vertices connected by $k$ parallel edges)
with $w_e = w$ for all $e$.
Then $Z_G(q,w) = q [q + (1+w)^k - 1]$,
so $Q_{\rm max}(G,\w) = |(1+w)^k - 1|$.
Now, if $|1+w| \ge 1$ we have
$\widehat\Delta(G,\w) = k |w||1+w|^{{k\over 2}-1}$ and $\Psi(G,\w) = |1+w|^k$.
Therefore, as $|w| \to \infty$ at fixed $k$,
the bound of Theorem~\ref{thm1.2c}
is a factor $4k$ from being sharp.

On the other hand, we may first apply parallel reduction
to yield a simple graph $\widehat{G} = K_2$
with weight $\widehat{w} = (1+w)^k - 1$ on its single edge,
and then apply Theorem~\ref{thm1.2c} or \ref{thm1.3b}
to $(\widehat{G},\widehat{\w})$.
The resulting bound is then (as $|w| \to \infty$)
a factor 4 or $\approx 4.892888$ from being sharp
(see Example~\ref{exam.K2}).
\qed
\end{example}

\begin{example}
   \label{exam.Cn}
\rm
Let $G$ be the $n$-cycle $C_n$ (which is simple for $n \ge 3$),
with $w_e = w$ for all $e$.
Then $Z_G(q,w) = (q+w)^n + (q-1) w^n$.
As $|w| \to \infty$ at fixed $n$,
we have $Q_{\rm max}(G,\w) = |w|^{n/(n-1)} + O(|w|)$.
On the other hand, if $|1+w| \ge 1$ we have
$\widehat\Delta = 2|w|$, $\Delta^* = 2|w| \, |1+w|^{1/2}$
and $\Psi = |1+w|^2$.
Therefore, as $|w| \to \infty$ the bounds of
Theorems~\ref{thm1.2c} and \ref{thm1.3b}
are $8|w|^{2} + O(|w|)$ and $2K^*_0 |w|^{3/2} + O(|w|)$, respectively
(here $2K^*_0  \approx 9.785776$).
Both of these bounds have the wrong order of magnitude
as $|w| \to \infty$ at fixed $n \ge 4$,
but the bound given by Theorem~\ref{thm1.3b} is a significant
improvement over that given by Theorem~\ref{thm1.2c}.
\qed
\end{example}


\begin{example}
  \label{exam.Kn}
\rm
Let $G$ be the complete graph $K_n$.
Take $w_e = w > 0$ for all $e$, with $w$ fixed independent of $n$
(unlike the usual \cite{Bollobas_96} scaling $w = \lambda/n$).
Then Janson \cite{Janson_08} has very recently proven that
\be
   \lim_{n \to \infty}
   {1 \over n^2} \log Z_{K_n}(e^{\alpha n},w)
   \;=\;
   \max[\smhalf \log(1+w), \alpha]
   \qquad\hbox{for } \alpha \ge 0
   \;.
 \label{eq.janson}
\ee
[This is because the sum \reff{def.ZG} is dominated by two contributions:
 the terms with $(V,A)$ connected, which together contribute
 $e^{\alpha n} (1+w)^{\binom{n}{2}} [1+o(1)]$,
 and the term $A = \emptyset$, which contributes $e^{\alpha n^2}$.]
It then follows from the Yang--Lee \cite{Yang-Lee_52}
theory of phase transitions
(see e.g.\ \cite[Theorem~3.1]{Sokal_chromatic_roots})
that $Z_{K_n}(e^{\alpha n},w)$ must have complex roots $\alpha_n$
that converge to $\alpha_\star = \smhalf \log(1+w)$
as $n \to\infty$.
Hence $Q_{\rm max}(K_n,\w) \ge (1+w)^{n/2 + o(n)}$
(and this is presumably the actual order of magnitude).
On the other hand, we have $\Delta^*(K_n,\w) = (n-1)w(1+w)^{n/2-1}$,
so that the upper bound given by Theorem~\ref{thm1.3b} is nearly sharp
when $n \to\infty$ at fixed $w > 0$
[it exceeds the truth by at most a factor $e^{o(n)}$
 even though both the truth and the bound are growing exponentially in $n$].

By contrast, $\widehat\Delta=(n-1)w(1+w)^{(n-3)/2}$ and $\Psi=(1+w)^{n-1}$,
so the bound of Theorem~\ref{thm1.2c} is much worse
because of its growth as $(1+w)^{n-2}$ rather than $(1+w)^{n/2-1}$.
\qed
\end{example}

\begin{example}
  \label{exam.Zd}
\rm
Let $G$ be a large finite piece of the simple hypercubic lattice $\Z^{d}$
(for some fixed $d \ge 2$) with nearest-neighbor edges,
and take $w_e = w > 0$ for all $e$.
For real $q > 0$ sufficiently large,
it is known \cite{Martirosian_86,Laanait_86,Kotecky_90,Laanait_91,Borgs_91}
that the first-order phase-transition point $w_t$ lies at
\be
   w_t(q)   \;=\;  q^{1/d} \,+\, O(1)  \;.
\ee
It then follows from the Yang--Lee \cite{Yang-Lee_52}
theory of phase transitions
that there will be complex zeros of the partition function
arbitrarily close (as $G$ grows) to the phase-transition point $(q,w_t(q))$;
so as $w \uparrow \infty$ (for fixed $d \ge 2$) we will have asymptotically
$Q_{\rm max}(G,\w) \ge w^d [1 + O(1/w)]$
(and this is presumably the actual order of magnitude).
Since
$\Delta^*(G,\w) = 2dw(1+w)^{d-1/2}$,
the upper bound given by Theorem~\ref{thm1.3b} is
off by at most a factor of order $w^{1/2}$
(i.e.\ it grows as $w^{d+{1/2}}$ instead of $w^d$).
By contrast, $\widehat\Delta=2dw(1+w)^{d-1}$ and $\Psi=(1+w)^{2d}$, so the bound of Theorem~\ref{thm1.2c} is again much worse,
because it grows as $w^{2d}$ rather than $w^{d+1/2}$.
\qed
\end{example}

\begin{example}
  \label{exam.G1G2}
\rm
Let $G$ be a disjoint union $G = G_1 \uplus G_2$.
Then $Q_{\rm max}(G) = \max\{ Q_{\rm max}(G_1),$ $Q_{\rm max}(G_2) \}$,
$\widehat{\Delta}(G) = \max\{ \widehat{\Delta}(G_1), \widehat{\Delta}(G_2) \}$
and $\Psi(G) = \max\{ \Psi(G_1), \Psi(G_2) \}$.
But the product $\widehat\scrk(\Psi) \, \widehat\Delta$ for $G$
{\em can exceed}\/ the maximum of those for $G_1$ and $G_2$
because one factor could be maximized for $G_1$ and the other for $G_2$.
For instance, for $i=1,2$ let $G_i$ be an $r_i$-regular graph
with all edge weights equal to $w_i$, where $|1+w_i| \ge 1$.
Then
\begin{subeqnarray}
   \widehat{\Delta}(G_i)  & = &  r_i \, |w_i| \, |1+w_i|^{r_i/2 - 1}  \\[2mm]
   \Psi(G_i)     & = &  |1+w_i|^{r_i}
\end{subeqnarray}
Now choose (for instance)
$r_1 = \rho \gg 1$, $r_2 = 3$, $w_1 = 1$, $w_2 \gg 1$.
Then
\be
   {\widehat{\Delta}(G_1) \over \widehat{\Delta}(G_2)}
   \;=\;
   {\rho \, 2^{\rho/2 - 1}  \over  3 \, w_2 (1+w_2)^{1/2}}
   \;\approx\;
   {\rho \, 2^{\rho/2}  \over  6 \, w_2^{3/2}}
\ee
while
\be
   {\Psi(G_2) \over \Psi(G_1)}
   \;=\;
   {(1+w_2)^3  \over 2^\rho}
   \;\approx\;
   {w_2^3  \over 2^\rho}
   \;.
\ee
So if we choose
\be
   \rho^2 2^\rho  \;\gg\;  w_2^3  \;\gg\;  2^\rho
\ee
we will have $\widehat{\Delta}(G_1) \gg \widehat{\Delta}(G_2)$
but $\Psi(G_2) \gg \Psi(G_1)$.
\qed
\end{example}


\appendix
\section{Appendix: Proof of Lemma~\ref{lemma.oldcor.borgs} and related facts}
\zeq

In this appendix we prove Lemma~\ref{lemma.oldcor.borgs}.
Actually, we prove much more:
though only parts (e,f,h) of Proposition~\ref{prop.borgs} below
actually arise in Lemma~\ref{lemma.oldcor.borgs}
and hence in the proofs of Theorems~\ref{thm1.2c} and \ref{thm1.3b},
we think it worthwhile to collect here some
additional properties of the function $F_\lambda(\beta)$
defined by \reff{prop.borgs.def_F}.
Some of these properties will be invoked in the Discussion
after the proof of Theorem~\ref{thm1.3b},
while others may end up playing a role in future work.

\begin{proposition}
  \label{prop.borgs}
For $\lambda \ge 0$ and $\beta > 0$, define the function
\be
   F_\lambda(\beta)   \;=\;
   \min\left\{ L  \colon\;
                  \inf\limits_{\alpha > 0} \,
   (e^\alpha-1)^{-1} \sum\limits_{n=2}^\infty  e^{\alpha n} \,
         L^{-(n-1)} \, {[1 + (n-1)\lambda]^{n-2} \over (n-1)!}
   \;\le\;   \beta
       \right\} \;. \qquad
 \label{prop.borgs.def_F}
\ee
Then:
\begin{itemize}
   \item[(a)]  $F_\lambda(\beta)$ is an increasing function of $\lambda$
and a decreasing function of $\beta$.
   \item[(b)]  $\beta F_\lambda(\beta)$ is an increasing function
of both $\lambda$ and $\beta$.
   \item[(c)]  $F_\lambda(\mu/\lambda)/\lambda$
is a decreasing function of both $\lambda$ and $\mu \, (> 0)$.
In particular, $F_\lambda(\beta)/\lambda$
is a decreasing function of both $\lambda$ and $\beta$.
   \item[(d)] $\log F_\lambda(\beta)$ is a convex function of $\log\beta$.
   \item[(e)]  We have
\be
   F_\lambda(\beta)
   \;=\;
   \min\limits_{1 < y < 1+\beta} \,
        {\beta y^{\lambda} \over (1+\beta-y) \log y}
   \;.
  \label{prop.borgs.F_a}
\ee
   \item[(f)] For $\lambda=0,1$ we have
\begin{eqnarray}
   F_0(\beta)
   & = &
     {\beta \over 1+\beta} \,
     W( (1+\beta)e )
     \big/
     \left[ W( (1+\beta)e ) \,-\, 1 \right]^2
      \label{f0}  \\[2mm]
%
   F_1(\beta)
   & = &
     \beta\,W\biggl( \displaystyle {e \over 1+\beta} \biggr)
     \bigg/
      \left[ 1 - W\biggl( \displaystyle {e \over 1+\beta}
                  \biggr) \right]^2
 \label{f1}
\end{eqnarray}
where $W$ is the Lambert $W$ function \cite{Corless_96},
 i.e.\ the inverse function to $x \mapsto x e^x$.
   \item[(g)] For $0 \le \lambda \le \lambda'$ we have
\be
   F_\lambda(\beta)  \;\le\;
     {1+2\lambda\hphantom{{}'} \over 1+2\lambda'} \;
     F_{\lambda'} \!\biggl( {1+2\lambda\hphantom{{}'} \over 1+2\lambda'}
                            \, \beta \biggr)
   \;.
 \label{eq.Flambda.lambdaprime.bound}
\ee
   \item[(h)] For $0 \le \lambda \le 1$ we have
\be
   F_\lambda(\beta)  \;\le\;  4\beta^{-1} + (1+2\lambda)
   \;.
 \label{eq.Flambda.bound}
\ee
\end{itemize}
\end{proposition}

\medskip

\proofof{Proposition~\ref{prop.borgs}}
(a) It is immediate from the definition \reff{prop.borgs.def_F}
that $F_\lambda(\beta)$ is increasing in $\lambda$ and decreasing in $\beta$.

(b) The change of variables $L'= \beta L$
in \reff{prop.borgs.def_F} shows that
\be
   \beta F_\lambda(\beta)
   \;=\;
   \min\left\{ L'  \colon\;
                  \inf\limits_{\alpha > 0} \,
   (e^\alpha-1)^{-1} \sum\limits_{n=2}^\infty  e^{\alpha n} \,
         (L')^{-(n-1)} \, \beta^{n-2} \,
         {[1 + (n-1)\lambda]^{n-2} \over (n-1)!}
   \;\le\;   1
       \right\}
\ee
is increasing in both $\lambda$ and $\beta$.

(c) The change of variables $L''=L/\lambda$
in \reff{prop.borgs.def_F} shows that
\be
   {F_\lambda(\mu/\lambda)  \over  \lambda}
   \;=\;
   \min\left\{ L''  \colon\;
                  \inf\limits_{\alpha > 0} \,
   (e^\alpha-1)^{-1} \sum\limits_{n=2}^\infty  e^{\alpha n} \,
         (L'')^{-(n-1)} \, {[\lambda^{-1} + (n-1)]^{n-2} \over (n-1)!}
   \;\le\;   \mu
       \right\}
\ee
is decreasing in both $\lambda$ and $\mu$.

(d) Suppose that we have triplets $(\alpha_i, L_i, \beta_i)$
satisfying
\be
   \sum_{n=2}^\infty  e^{\alpha_i n} \, L_i^{-(n-1)} \,
      {[1+(n-1)\lambda]^{n-2} \over (n-1)!}
   \;\le\;
   \beta_i \, (e^{\alpha_i} - 1)
\ee
for $i=1,2$.  Now let $\kappa \in [0,1]$ and define
\begin{subeqnarray}
   \bar{\alpha}  & = &  \kappa \alpha_1 + (1-\kappa) \alpha_2   \\
   \bar{L}       & = &  L_1^\kappa  L_2^{1-\kappa}              \\
   \bar{\beta}   & = &  \beta_1^\kappa  \beta_2^{1-\kappa}
\end{subeqnarray}
Then H\"older's inequality with $p=1/\kappa$ and $q=1/(1-\kappa)$ yields
\be
   \sum_{n=2}^\infty  e^{\bar{\alpha} n} \: \bar{L}^{-(n-1)} \:
      {[1+(n-1)\lambda]^{n-2} \over (n-1)!}
   \;\le\;
   \bar{\beta} \, (e^{\alpha_1} - 1)^\kappa \, (e^{\alpha_2} - 1)^{1-\kappa}
   \;.
\ee
And since the function $\alpha \mapsto \log(e^{\alpha} - 1)$
is concave on $(0,\infty)$, we have
$(e^{\alpha_1} - 1)^\kappa  (e^{\alpha_2} - 1)^{1-\kappa}
  \le e^{\bar{\alpha}} - 1$.
This proves (d).

(e) The proof that \reff{prop.borgs.def_F} is equivalent to
\reff{prop.borgs.F_a} will be modelled on an
argument of Borgs \cite[eq.~(4.22)~ff.]{Borgs_06},
who proved a related result.

Note first that $c \mapsto ce^{-c}$ maps the interval $[0,1]$
strictly monotonically onto the interval $[0,1/e]$;
and recall \cite[p.~28]{Stanley_99} that its inverse map
is the tree function
\be
   T(x)  \;=\;  \sum_{n=1}^\infty {n^{n-1} \over n!} \, x^n
   \;,
\ee
which is convergent and monotonically increasing for $0 \le x \le 1/e$
and satisfies $T(c e^{-c}) = c$ for $0 \le c \le 1$.
Moreover, it is well known (see e.g. \cite[eq.~(2.36)]{Corless_96})
that for all real $\kappa>0$ one has [cf.\ \reff{gen_fn}]
\be
   \left({T(z)\over z}\right)^\kappa
   \;=\;
   \sum_{m=0}^\infty {\kappa\,(m+\kappa)^{m-1} \over m!} \, z^m
\ee
(this is an easy consequence of the Lagrange inversion formula).
Writing for convenience $U(z) = T(z)/z$, we therefore have
\be
   \sum_{n=1}^\infty {[1+ (n-1)\lambda]^{n-2} \over (n-1)!} \, z^n
   \;=\;
   z \, U(\lambda z)^{1/\lambda}
\ee
for all real $\lambda > 0$.

The inequality on the right-hand side of \reff{prop.borgs.def_F}
is then equivalent to the statement that $\lambda e^\alpha/L \le 1/e$
(otherwise the sum would be divergent) and
\be
   e^\alpha \, U(\lambda e^{\alpha}/L)^{1/\lambda} \,-\, e^\alpha
   \;\le\;  \beta (e^\alpha - 1)
   \;.
 \label{eq.borgs.1}
\ee
Eliminating $L$ in favor of a new variable $c$ defined by
$\lambda e^{\alpha}/L=ce^{-c}$ with $0 \le c \le 1$,
and using the fact that $U(c e^{-c}) = e^c$,
we see that
the inequality on the right-hand side of \reff{prop.borgs.def_F}
is equivalent to
\be
   c \;\le\;  \min\Bigl\{ 1,\, \lambda \log[1 + \beta (1-e^{-\alpha})] \Bigr\}
   \;.
 \label{ineq.c}
\ee
Since $L = \lambda e^{\alpha}/(ce^{-c})$,
and $ce^{-c}$ increases monotonically with $c$ for $0 \le c \le 1$,
we deduce that \reff{ineq.c}
is equivalent to
\be
   L  \;\ge\;
   \begin{cases}
      \displaystyle
      {e^\alpha \, [1 + \beta (1-e^{-\alpha})]^\lambda  \over
       \log[1 + \beta (1-e^{-\alpha})]}
           & \text{if } \beta (1-e^{-\alpha}) \le e^{1/\lambda} -1  \\[5mm]
      \lambda e^{\alpha+1}
           & \text{if } \beta (1-e^{-\alpha}) \ge e^{1/\lambda} -1
   \end{cases}
 \label{ineq.L.1}
\ee
Changing variables from $\alpha$ to
$y = 1 + \beta (1-e^{-\alpha})$,
we can rewrite this as
\be
   L  \;\ge\;
   \begin{cases}
      \displaystyle
      {\beta y^{\lambda} \over (1+\beta-y) \log y}
            & \text{if } 1 < y < \min(e^{1/\lambda}, 1+\beta)  \\[6mm]
      \displaystyle
      {\lambda\beta e \over 1+\beta-y}
            & \text{if } e^{1/\lambda} \le y < 1+\beta
   \end{cases}
 \label{ineq.L.2}
\ee
Now we can optimize over $y$:
the minimum will always be found in the interval $1 < y \le e^{1/\lambda}$,
so we have
\be
   F_\lambda(\beta)
   \;=\;
   \min\limits_{1 < y < \min(e^{1/\lambda}, 1+\beta)} \,
        {\beta y^{\lambda} \over (1+\beta-y) \log y}
   \;=\;
   \min\limits_{1 < y <  1+\beta} \,
        {\beta y^{\lambda} \over (1+\beta-y) \log y}
   \;,
\ee
where the final equality results from the fact that
$y^{\lambda}/[(1+\beta-y) \log y]$ is increasing for
$e^{1/\lambda} \le y < 1+\beta$.
This proves the equivalence of \reff{prop.borgs.def_F}
with \reff{prop.borgs.F_a} for $\lambda > 0$;
and the case $\lambda=0$ follows by taking limits
(or by an easy direct proof).

(f) For $\lambda=0$, simple calculus shows that
the minimum in \reff{prop.borgs.F_a} is attained at
$y = (1+\beta) / W((1+\beta)e)$,
so that $F_0(\beta)$ is given by \reff{f0}.
Likewise, for $\lambda=1$, simple calculus shows that
the minimum in \reff{prop.borgs.F_a} is attained at
$y = (1+\beta) \, W(e/(1+\beta))$,
so that $F_1(\beta)$ is given by \reff{f1}.

(g) To prove the comparison inequality \reff{eq.Flambda.lambdaprime.bound},
it suffices to observe that whenever $0 \le \lambda \le \lambda'$
and $n \ge 2$ we have
\be
   \biggl( {1+(n-1)\lambda\hphantom{{}'}  \over 1+(n-1)\lambda'} \biggr)^{n-2}
   \;\le\;
   \biggl( {1+2\lambda\hphantom{{}'}  \over 1+2\lambda'} \biggr)^{n-2}
\ee
(just consider $n=2$ and $n \ge 3$ separately).
Inserting this into the definition \reff{prop.borgs.def_F}
yields \reff{eq.Flambda.lambdaprime.bound}.

(h) To prove the upper bound \reff{eq.Flambda.bound},
write $y=1+x$ in \reff{prop.borgs.F_a} and use the inequalities
\begin{eqnarray}
   {1 \over \log(1+x)}  & \le &  {1 \over x} \,+\, {1 \over 2}
      \label{eq.fbound} \\[3mm]
   (1+x)^\lambda   & \le &  1 + \lambda x
\end{eqnarray}
which are valid for all $x > 0$ and $0 \le \lambda \le 1$.\footnote{
   {\sc Proof of \reff{eq.fbound}:}
   Write $t = \log(1+x) > 0$;  then \reff{eq.fbound} states that
   $1/t \le 1/(e^t-1) + 1/2$.
   This is trivially true for $t \ge 2$;
   and for $0 < t < 2$ it is equivalent to $e^t - 1 \le t/(1-t/2)$,
   which is obvious from the Taylor series.
}
Therefore,
\be
   {\beta y^\lambda \over (1+\beta-y) \log y}
   \;\le\;
   {\beta (1+ \lambda x) \bigl( {1 \over x} \,+\, {1 \over 2} \bigr)
    \over  \beta - x}
   \;.
\ee
The latter function is minimized at
$x = (-2+ \sqrt{4+(2+4\lambda)\beta+2\lambda\beta^2})/[1+(2+\beta)\lambda]
 \in (0,\beta)$,
with minimum value
\be
   {1 \over 2} \,+\, \lambda \,+\, {2 \over \beta}
     \,+\, {2 \over \beta} \sqrt{(1+\beta/2)(1+\lambda\beta)}
   \;.
\ee
This, in turn, is bounded above by $4\beta^{-1} + (1+2\lambda)$
on the entire interval $0 < \beta < \infty$.\footnote{
   {\sc Proof:}
   We have
   $$ \sqrt{(1+c_1 \beta)(1+c_2 \beta)}  \;\le\;
      1 \,+\, {c_1 + c_2 \over 2} \beta
   $$
   for all $c_1,c_2,\beta \ge 0$,
   as is easily seen by squaring both sides and using
   the arithmetic-geometric-mean inequality
   $\sqrt{c_1 c_2} \le (c_1 + c_2)/2$.
}
[Alternatively, it suffices to make this proof for $\lambda=1$
 and then invoke \reff{eq.Flambda.lambdaprime.bound}
 to deduce the result for $0 \le \lambda < 1$.]
\qed


\bigskip

{\bf Remarks.}
1.  The proof of Proposition~\ref{prop.borgs}(e)
becomes a bit simpler for $\beta \le e^{1/\lambda}-1$,
since we then always have $\beta (1-e^{-\alpha}) \le e^{1/\lambda}-1$
and hence we need not worry about the second case in
\reff{ineq.L.1} and \reff{ineq.L.2}.
This simplification applies in particular when
$\lambda \le 1$ and $\beta \le 1$,
which covers what is needed in the proofs of both
Theorem~\ref{thm1.2c} ($\lambda=1$, $\beta = \psi^{-1/2} \le 1$)
and Theorem~\ref{thm1.3b} ($0 < \lambda \le 1$, $\beta = 1$).

2.  We can compute the small-$\beta$ asymptotics of $F_\lambda(\beta)$
by expanding \reff{prop.borgs.F_a} in powers of $y-1$:
the minimum is located at
\be
   y \;=\;  1  \,+\, {1 \over 2} \beta  \,-\, {1+2\lambda \over 16} \beta^2
               \,+\, {5+12\lambda \over 192} \beta^3
               \,-\, {43+122\lambda+12\lambda^2-24\lambda^3 \over 3072} \beta^4
               \,+\, \ldots
   \quad
\ee
and we have
\be
   F_\lambda(\beta)  \;=\;
   4\beta^{-1} \,+\, (1+2\lambda)
               \,-\, {7+12\lambda-12\lambda^2 \over 48} \beta
               \,+\, {11+26\lambda-12\lambda^2-8\lambda^3 \over 192} \beta^2
               \,+\, \ldots
   \;.
\ee
For $\lambda=0,1$ an alternate method is to expand \reff{f0}/\reff{f1}:
we obtain
\begin{eqnarray}
   F_0(\beta)
   & = &
   4\beta^{-1} \,+\, 1 \,-\, \frac{7}{48} \beta \,+\, \frac{11}{192} \beta^2
      \,-\, \frac{443}{15360} \beta^3 \,+\, \frac{607}{36864} \beta^4
      \,-\, \ldots
      \quad \\[3mm]
   F_1(\beta)
   & = &
   4\beta^{-1} \,+\, 3 \,-\, \frac{7}{48} \beta \,+\, \frac{17}{192} \beta^2
      \,-\, \frac{923}{15360} \beta^3 \,+\, \frac{8113}{184320} \beta^4
      \,-\, \ldots
   \quad
\end{eqnarray}
Therefore, the large-$\psi$ asymptotics of
$\widehat\scrk(\psi) = F_1(\psi^{-1/2})$ is
\be
   \widehat\scrk(\psi)  \;=\;
   4\psi^{1/2} \,+\, 3 \,-\, \frac{7}{48} \psi^{-1/2}
      \,+\, \frac{17}{192} \psi^{-1} \,-\, \frac{923}{15360} \psi^{-3/2}
      \,+\, \frac{8113}{184320} \psi^{-2}
      \,-\, \ldots
  \;. \quad
 \label{eq.Kpsi.largepsi}
\ee

3.  In the preprint version of this paper\footnote{
   \url{http://arxiv.org/abs/0810.4703v2}
},
we conjectured (based on plots of $F_1$ and its derivatives)
that $F_1(\beta)$ is a {\em completely monotone}\/ function of $\beta$
on $(0,\infty)$,
i.e.\ $(-1)^k \, d^{\,k}\!F_1(\beta)/d\beta^k \ge 0$ for all $\beta > 0$
and all integers $k \ge 0$,
and indeed that $G_1(\beta) = F_1(\beta) - 4/\beta$
is completely monotone, which is stronger.\footnote{
   See e.g.\ \cite{Widder_46} for the theory of completely monotone
   functions on $(0,\infty)$ --- in particular the
   Bernstein--Hausdorff--Widder theorem,
   which states that a function is completely monotone on $(0,\infty)$
   if and only if it is the Laplace transform of a positive measure
   supported on $[0,\infty)$.
}
Even more strongly, we conjectured
(based on computations for $\imag\beta > 0$)
that $G_\lambda(\beta) = F_\lambda(\beta) - 4/\beta$
is a {\em Stieltjes function}\/ for $\lambda=0$ and $\lambda=1$,
i.e.\ it can be written in the form
\be
   f(\beta)  \;=\;  C \,+ \int\limits_{[0,\infty)} \! {d\rho(t) \over \beta+t}
 \label{def.stieltjes}
\ee
where $C \ge 0$ and $\rho$ is a positive measure on $[0,\infty)$.\footnote{
   More information on Stieltjes functions can be found in
   \cite{Widder_46} \cite[pp.~126--128]{Akhiezer_65}
   \cite{Berg_79,Berg_80,Berg_08,Sokal_stieltjes}
   and the references cited therein.
   In order to test numerically the Stieltjes property
   for $G_0(\beta)$ and $G_1(\beta)$,
   we have used the complex-analysis characterization:
   a function $f \colon\, (0,\infty) \to \R$
   is Stieltjes if and only if it is the restriction to $(0,\infty)$
   of an analytic function on the cut plane $\C \setminus (-\infty,0]$
   satisfying $f(z) \ge 0$ for $z>0$ and
   $\Im f(z) \le 0$ for $\Im z > 0$.
   See e.g.\ \cite[p.~127]{Akhiezer_65} or \cite{Berg_80}.
}
This latter conjecture has now been proven by Kalugin, Jeffrey and Corless
\cite{Kalugin_11}.
It is even possible that $G_\lambda$ is a Stieltjes function
also for $0 < \lambda < 1$, but a different method of proof will be needed.
\qed

\section*{Acknowledgments}

We are extremely grateful to Svante Janson for answering
our query about the behavior of $Z_{K_n}(q,w)$ when $w > 0$
is taken independent of $n$ [cf.\ \reff{eq.janson}].
We are also indebted to an anonymous referee,
whose incisive comments on a previous version of our paper
led us to obtain some notable improvements.

We wish to thank the Isaac Newton Institute for Mathematical Sciences,
University of Cambridge, for generous support during the programme on
Combinatorics and Statistical Mechanics (January--June 2008),
where this work was begun.
One of us (A.D.S.)\ also thanks the
Institut Henri Poincar\'e -- Centre Emile Borel
for hospitality during the programmes on
Interacting Particle Systems, Statistical Mechanics and Probability Theory
(September--December 2008)
and Statistical Physics, Combinatorics and Probability
(September--December 2009),
and the Laboratoire de Physique Th\'eorique
at the \'Ecole Normale Sup\'erieure
for hospitality during April--June 2011.

This research was supported in part by
U.S.\ National Science Foundation grant PHY--0424082,
by the Conselho Nacional de Desenvolvimento Cient\'{\i}fico e
Tecnol\'ogico (CNPq), and by the
Funda\c{c}\~ao de Amparo \`a Pesquisa do Estado de Minas Gerais (FAPEMIG).

\end{document}